\documentclass[12pt,draftcls,onecolumn]{IEEEtran}

\usepackage{amsfonts}
\usepackage{amssymb}
\usepackage{amsmath}
\usepackage{mathptmx}
\usepackage{cancel}
\usepackage{pstricks}
\usepackage{psfrag}
\usepackage{mathrsfs}
\usepackage{graphicx}
\usepackage{pgf,tikz}
\usepackage{algorithm}
\usepackage{algpseudocode}
\usepackage{soul}
\def\be{\begin{equation}}
\def\ee{\end{equation}}
\def\bea{\begin{eqnarray}}
\def\eea{\end{eqnarray}}
\def\beann{\begin{eqnarray*}}
\def\eeann{\end{eqnarray*}}

\newcommand{\rank}{{\rm rank}}
\newcommand{\nn}{{\nonumber}}

\def\ns{\hspace{-1mm}}
\newcommand{\real}{{\mathbb{R}}}

\def\spacingset#1{\def\baselinestretch{#1}\small\normalsize}
\spacingset{1.35}

\newtheorem{lemma}{Lemma}
\newtheorem{theorem}{Theorem}
\newtheorem{remark}{Remark}
\newtheorem{corollary}{Corollary}
\newtheorem{proposition}{Proposition}
\newtheorem{definition}{Definition}
\newtheorem{problem}{Problem}
\newtheorem{example}{Example}[section]
\newtheorem{assumption}{Assumption}[section]
\newenvironment{proofof}{{\em Proof of }}{\hfill \hspace*{1pt}\hfill $\blacksquare$}

\def\be{\begin{equation}}
\def\ee{\end{equation}}
\def\bea{\begin{eqnarray}}
\def\eea{\end{eqnarray}}
\def\beann{\begin{eqnarray*}}
\def\eeann{\end{eqnarray*}}

\def\ns{\hspace{-1mm}}

\def\proof{\noindent{\bf{\em Proof:}\ \ }}
\def\QED{\mbox{\rule[0pt]{1.5ex}{1.5ex}}}
\def\endproof{\hspace*{\fill}~\QED\par\endtrivlist\unskip}

\newcommand{\ima}{\operatorname{im}}
\newcommand{\normrank}{\operatorname{normrank}}
\newcommand{\diag}{\operatorname{diag}}
\newcommand{\defi}{\stackrel{\text{\tiny def}}{=}}

\newcommand{\complex}{{\mathbb{C}}}

\def\gA{{\cal A}}
\def\gB{{\cal B}}

\def\gD{{\cal D}}
\def\gE{{\cal E}}

\def\gG{{\cal G}}
\def\gH{{\cal H}}

\def\gJ{{\cal J}}
\def\gK{{\cal K}}
\def\gL{{\cal L}}
\def\gM{{\cal M}}

\def\gP{{\cal P}}

\def\gR{{\cal R}}
\def\gS{{\cal S}}
\def\gT{{\cal T}}
\def\gU{{\cal U}}
\def\gV{{\cal V}}
\def\gW{{\cal W}}
\def\gX{{\cal X}}
\def\gY{{\cal Y}}
\def\gZ{{\cal Z}}

\def\bmat{\left[ \begin{array}}
\def\emat{\end{array} \right]}

\def\bmat{\left[ \begin{array}}
\def\emat{\end{array} \right]}
\def\bsmat{\left[ \begin{smallmatrix}}
\def\esmat{\end{smallmatrix} \right]}

\def\l{{\lambda}}
\def\gA{{\cal A}}
\def\gP{{\cal P}}
\def\gB{{\cal B}}
\def\gU{{\cal U}}
\def\gL{{\cal L}}
\def\gR{{\cal R}}
\def\gG{{\cal G}}
\def\gV{{\cal V}}
\def\gH{{\cal H}}
\def\gS{{\cal S}}
\def\gK{{\cal K}}
\def\gT{{\cal T}}

\def\gX{{\cal X}}
\def\i{{i}}
\newcommand{\spanR}{\operatorname{span}}

\begin{document}
\title{\LARGE{Globally monotonic tracking control\\
 of multivariable systems}}

\author{Lorenzo~Ntogramatzidis, Jean-Fran\c{c}ois Tr{\'e}gou{\"e}t and Robert Schmid 
 \thanks{L. Ntogramatzidis is with the Department of Mathematics and
Statistics, Curtin University, Perth,
Australia. E-mail: {\tt L.Ntogramatzidis@curtin.edu.au. }}

 \thanks{Jean-Fran\c{c}ois Tr{\'e}gou{\"e}t is with Universit\'{e} de Lyon, Laboratoire Amp\`{e}re CNRS UMR 5005, INSA-Lyon; F-69621, Villeurbanne, France. E-mail: {\tt jean-francois.tregouet@insa-lyon.fr. } (Research partially carried out at Curtin University).}

\thanks{R. Schmid is with the Department of Electrical and Electronic Engineering,
The University of Melbourne, Parkville, VIC 3010, Australia. E-mail: {\tt rschmid@unimelb.edu.au. }}

\thanks{Partially supported by the Australian Research Council under the grant FT120100604.}
}

\markboth{DRAFT}{Shell \MakeLowercase{\textit{et al.}}: Bare Demo of IEEEtran.cls for Journals}

\maketitle

\vspace{-1cm}

\IEEEpeerreviewmaketitle

\begin{abstract}
%
In this paper we present a method for designing a linear time invariant (LTI) state-feedback controller to monotonically track a constant step reference at any desired rate of convergence for any arbitrarily assigned initial condition.
Necessary and sufficient constructive conditions are given to deliver a monotonic step response from all initial conditions. This method is developed for multi-input multi-output (MIMO) systems, and can be applied to square and non-square systems, strictly proper and bi-proper systems, and, importantly, also minimum and non-minimum phase systems. The control methods proposed here show that for MIMO LTI systems the objectives of achieving a rapid settling time, while at the same time
avoiding overshoot and/or undershoot, are not necessarily competing objectives.
\end{abstract}


\section{Introduction}
\label{secintro}
The problem of improving the shape of the step response curve for linear time invariant (LTI) systems is as old as control theory. Its relevance is seen in countless applications such as heating/cooling systems, elevator and satellite positioning, automobile cruise control and the positioning of a CD disk read/write head. The common element in these problems involves designing a control input for the system to make the output take a certain desired target value, and then keep it there.
 
 A fundamental issue in classical feedback control is the design of control laws that provide good performance both at steady state and during the transient. The \textit{steady state performance} is typically assumed to be satisfactory if, once the transient vanishes, the output of the system is constant and equal (or very close) to the desired value.
When dealing with the \textit{transient performance}, one is usually concerned with the task of reducing both the overshoot and the undershoot, or, ideally, of achieving a monotonic response that rapidly converges to the steady-state value.
 It is commonly understood that the objectives of achieving a rapid (short) settling time, while at the same time avoiding overshoot and undershoot, are competing objectives in the controller design, and must be dealt with by seeking a trade-off, see e.g. \cite{Franklin-PE-94,Dorf-B-08}, or any standard textbook on the topic.
While this is certainly the case for single-input single-output (SISO) systems, the control methods we develop and implement in this paper challenge this widely-held perception for the multi-input multi-output (MIMO) case. We show in particular that in the case of LTI MIMO systems, it is possible to achieve arbitrarily fast settling time and also a monotonic step response in all output components for any initial condition, which naturally imply the avoidance of overshoot/undershoot even in the presence of non-minimum phase invariant zeros. 

In contrast with the extensive literature for SISO systems, which includes -- but is far from being limited to -- \cite{Hoagg-B-07,Lin-F-97,Stewart-D-06,Darbha-B-03,Bement-J-04,Darbha-03,Darbha-B-02,Middleton-91,Lau-MB-03} and the references cited therein, to date there have been very few papers offering analysis or design methods for avoiding undershoot or overshoot in the step response of MIMO systems, see e.g. \cite{Johansson-02} and the references therein. 
A recent contribution offering design methods for MIMO systems is \cite{Schmid-N-09}, where a procedure is proposed for the design of a state-feedback controller to yield a non-overshooting step response for LTI MIMO systems. Importantly, this design method is applicable to non-minimum phase systems, does not assume that the system state is initially at rest, and can be applied to both continuous-time and discrete-time (and also proper or bi-proper) systems. Very recently it has been shown in \cite{Schmid-N-09-1} how the method can be adapted to obtain a non-undershooting step response. 
The key idea behind the approach in \cite{Schmid-N-09} and \cite{Schmid-N-09-1} is 
to design the feedback matrix that achieves the desired closed-loop eigenstructure in such a way that only a small number of the closed-loop system modes appear in each component of the tracking error (which is defined as the difference between the system output and the desired target value). Indeed, if the closed-loop eigenstructure can be constrained in such a way that each component of the tracking error is driven only by a single real-valued closed-loop mode -- which is an exponential in the form $e^{\lambda\,t}$ in the continuous time or a power term $\lambda^k$ in the discrete time -- the output of the system is monotonic in each output component regardless of the initial condition of the system, and hence both overshoot and undershoot are avoided. For systems where the closed-loop eigenstructure can be constrained so that the error involves only the sum of two or three exponential terms (or powers in the discrete case) in each component, the design method offers a search algorithm for the selection of suitable closed-loop modes that ensures that the step response is non-overshooting, non-undershooting, or monotonic from any given initial condition and target reference. The feedback matrix is computed by inversion of the matrix of closed-loop eigenvectors.

However, a key limitations of the design methods given in \cite{Schmid-N-09} and \cite{Schmid-N-09-1} is the lack of analytic conditions, expressed in terms of the system structure, that guarantee the existence of a state-feedback controller that can deliver the desired transient response. More specifically, as already mentioned, the feedback matrix is computed by inversion of the matrix of closed-loop eigenvectors, and the problem is solvable if and only if this matrix is invertible, i.e., if and only if the set of closed-loop eigenvectors is linearly independent. When this is not the case, one may decide to change these eigenvectors by varying the choice of the corresponding closed-loop eigenvalues, and check whether the matrix of closed-loop eigenvectors has become non-singular as a result. However, it can very well happen that for any value of such eigenvalues, the corresponding eigenvectors are always linearly dependent. The method of \cite{Schmid-N-09} and \cite{Schmid-N-09-1} does not provide a structural criterion to decide if the problem admits a solution in terms of the problem data, nor does it guarantee that when the matrix of eigenvectors is singular one is allowed to conclude that the problem of achieving a monotonic response from any initial condition cannot be solved.
Moreover, as aforementioned, the design method involves a search for suitable closed-loop eigenvalues, and while this search can be conducted very efficiently, the authors were unable to give any conditions guaranteeing a satisfactory search outcome. 


 The objective of this paper is to revisit the design method of \cite{Schmid-N-09} and \cite{Schmid-N-09-1} to the end of developing
 conditions expressed in terms of the 
 system structure that are necessary and sufficient to guarantee that the design method will obtain a state-feedback controller that yields a monotonic step response from any initial condition and for any constant reference signal. When this goal is achievable, following \cite{Schmid-N-09}, we say that the control yields a {\em globally monotonic} response, by which we mean that the same state-feedback gain matrix yields a monotonic response from all initial conditions, and with respect to all possible step references. 
 
Thus, in this paper, for the first time in the literature, a complete and exhaustive answer to the problem of achieving a globally monotonic step response for a MIMO LTI system is provided. We show that, surprisingly, for MIMO LTI systems the presence of non-minimum phase invariant zeros does not prevent a globally monotonic step response to be achievable. Indeed, even in the presence of one or more non-minimum phase invariant zeros, it may still be possible to achieve a monotonic step response from any initial condition and for any constant reference signal. \\[-3mm]


{\bf Notation.} In this paper, the symbol $\{0\}$ stands for the origin of a vector space. 
For convenience, a linear mapping between finite-dimensional spaces and a matrix representation with respect to a particular basis are not distinguished notationally. The image and the kernel of matrix $A$ are denoted by $\ima\,A$ and $\ker\,A$, respectively. The Moore-Penrose pseudo-inverse of $A$ is denoted by $A^\dagger$. When $A$ is square, we denote by $\sigma(A)$ the spectrum of $A$. 
If $\gJ \subseteq \gX$, the restriction of the map $A$ to $\gJ$ is denoted by $A\,|\gJ$. If $\gX=\gY$ and $\gJ$ is $A$-invariant, the eigenvalues of $A$ restricted to $\gJ$ are denoted by $\sigma\,(A\,| \gJ)$. If $\gJ_1$ and $\gJ_2$ are $A$-invariant subspaces and $\gJ_1\,{\subseteq}\,\gJ_2$, the mapping induced by $A$ on the quotient space $\gJ_2 / \gJ_1$ is denoted by $A\,| {\gJ_2}/{\gJ_1}$, and its spectrum is denoted by $\sigma\,(A\,| {\gJ_2}/{\gJ_1})$.
The symbol $\oplus$ stands for the direct sum of subspaces. 
 The symbol $\uplus$ denotes union with any common elements repeated.
Given a map $A: \gX \longrightarrow \gX$ and a subspace $\gB$ of $\gX$, we denote by $\langle A, \gB \rangle$ the smallest $A$-invariant subspace of $\gX$ containing $\gB$. The symbol $\i$ stands for the imaginary unit, i.e., $\i=\sqrt{-1}$. Given a complex matrix $M$, the symbol $\overline{M}$ denotes the conjugate transpose of $M$. Moreover, we denote by $M_i$ its $i$-th row and by $M^j$ its $j$-th column, respectively.
{Given a Borel subset $\gG$ of $\real^n$, we denote the Lebesgue measure of $\gG$ in $\real^n$
as ${m}_{\scriptscriptstyle \real^n}(\gG)$. } 
 Given a finite set $S$, the symbol $2^S$ denotes the power set of $S$, while ${\rm {card}}(S)$ stands for the cardinality of $S$.
The symbols $\neg$, $\land$ and $\lor$ denote, respectively, the logical {\em not}, {\em and} and {\em or}. { Given a set $V$ and a logical proposition $\mathfrak{p}$ that depends on an object $v \in V$, we use the symbol $\{v\in V \,|\; \mathfrak{p}\}$ to indicate the subset in $V$ of the elements $v$ for which $\mathfrak{p}$ is true.}

\section{Problem Formulation}
\label{PF}
In what follows, whether the underlying system evolves in continuous or discrete time makes only minor differences and, accordingly, the time index set of any signal is denoted by $\mathbb{T}$, on the understanding that this represents either $\real^+$ in the continuous time or $\mathbb{N}$ in the discrete time. The symbol $\complex_g$ denotes either the open left-half complex plane $\complex^-$ in the continuous time or the open unit disc $\complex^\circ$ in the discrete time. A matrix $M \in \real^{n \times n}$ is said to be {\em asymptotically stable} if $\sigma(M) \subset \complex_g$, i.e., if it is Hurwitz in the continuous time and Schur in the discrete case. Finally, we say that $\l \in \complex$ is stable if $\l \in \complex_g$.
Consider the LTI system $\Sigma$ governed by
 \bea
 \Sigma: \
 \left\{ \begin{array}{lcr}
 \gD\,x(t) \ns&\ns = \ns&\ns A\,x(t)+B\,u(t),\;\;\;\; x(0)=x_0, \label{syseq1}\hfill\cr
 y(t) \ns&\ns = \ns&\ns C\,x(t)+D\,u(t),\hfill \end{array} \right.
 \label{sys}
 \eea
where, for all $t \in \mathbb{T}$, $x(t) \in \gX=\real^n$ is the
state, $u(t) \in \gU=\real^m$ is the control input, $y(t) \in \gY=\real^p$ is the output, and $A$, $B$, $C$ and $D$ are appropriate
dimensional constant matrices. The operator
$\gD$ denotes either the time derivative in the continuous time,
i.e., $\gD\, x(t)=\dot{x}(t)$, or the unit time shift
 in the discrete time, i.e., $\gD\, x(t)=x(t+1)$. 
 We assume with no loss of generality that all the columns of $\left[ \begin{smallmatrix} B \\[1mm] D \end{smallmatrix}\right]$ and all the rows of
 $[\, C \;\; D \,]$ are linearly independent. \footnote{If $\left[ \begin{smallmatrix} B \\[1mm] D \end{smallmatrix}\right]$ has non-trivial kernel, a subspace $\gU_0$ of the input space exists that does not influence the state dynamics. By performing a suitable (orthogonal) change of basis in the input space, we may eliminate $\gU_0$ and obtain an equivalent system for which this condition is satisfied. Likewise, if $[\, C \;\; D \,]$ is not surjective, there are some outputs that result as linear combinations of the remaining ones, and these can be eliminated using a dual argument using a change of coordinates in the output space.}

We recall that the Rosenbrock system matrix is defined as the matrix pencil
\bea
\label{ros}
P_{\scriptscriptstyle \Sigma}(\lambda) \defi \bmat{cc} A-\lambda\,I_n & B \\ C & D \emat
\eea
in the indeterminate $\lambda \in \complex$, see e.g. \cite{Rosenbrock-70}. The invariant zeros of $\Sigma$ are the values of $\lambda \in \complex$ for which the rank of $P_{\scriptscriptstyle \Sigma}(\lambda)$ is strictly smaller than its normal rank\footnote{The normal rank of a rational matrix $M(\lambda)$ is defined as
$\normrank M(\lambda) \defi \displaystyle\max_{\lambda \in \complex} \rank M(\lambda)$. The rank of $M(\lambda)$ is equal to its normal rank for all but finitely many $\lambda \in \complex$.}. More precisely,
 the invariant zeros are the roots of the non-zero polynomials on the principal diagonal of the Smith form of $P_{\scriptscriptstyle \Sigma}(\lambda)$, see e.g. \cite{Aling-S-84}.
Given an invariant zero $\lambda=z \in \complex$, the rank deficiency of $P_{\scriptscriptstyle \Sigma}(\lambda)$ at the value $\lambda=z$ is the geometric multiplicity of the invariant zero $z$, and is equal to the number of elementary divisors (invariant polynomials) of $P_{\scriptscriptstyle \Sigma}(\lambda)$ associated with the complex frequency $\lambda=z$. The degree of the product of the elementary divisors of $P_{\scriptscriptstyle \Sigma}(\lambda)$ corresponding to the invariant zero $z$ is the algebraic multiplicity of $z$, see \cite{MacFarlane-K-76}. The set of invariant zeros of $\Sigma$ is denoted by $\gZ$, and the set of minimum-phase invariant zeros of $\Sigma$ is 
$\gZ_g\defi \gZ\cap \complex_g$.

We denote by $\gV^\star$ the largest output-nulling subspace of $\Sigma$, i.e., the largest subspace $\gV$ of $\gX$ for which a matrix $F\,{\in}\,\mathbb{R}^{m\,{\times}\,n}$ exists such that $(A+B\,F)\,\gV\subseteq \gV \subseteq \ker (C+D\,F)$. Any real matrix $F$ satisfying this inclusion is called a {\it friend \/} of $\gV$. We denote by $\mathfrak{F}(\gV)$ the set of friends of $\gV$.
 The symbol $\gR^\star$ denotes the so-called {\em output-nulling reachability subspace} on $\gV^\star$, and is the smallest $(A\,{+}\,B\,F)$-invariant subspace of $\gX$ containing $\gV^\star\,{\cap}\,B\,\ker\,D$, where $F\,{\in}\,\mathfrak{F}(\gV^\star)$.
 The closed-loop spectrum can be partitioned as
$\sigma(A+B\,F)=\sigma(A+B\,F\,|\,\gV^\star)\uplus \sigma(A+B\,F\,|\,\gX/\gV^\star)$, where
 $\sigma(A+B\,F\,|\,\gV^\star)$ is the spectrum of $A+B\,F$ restricted to $\gV^\star$, and its elements are the {\em inner eigenvalues} of the closed-loop with respect to $\gV^\star$, and
$\sigma(A+B\,F\,|\,\gX/\gV^\star)$ is the spectrum of the mapping induced by $A+B\,F$ on the quotient space $\gX/\gV^\star$. Its elements are referred to as {\em outer eigenvalues} of the closed-loop with respect to $\gV^\star$. We say that $\gV^\star$ is {\em inner stabilisable} if a friend $F$ of $\gV^\star$ exists such that $\sigma(A+B\,F\,|\,\gV^\star) \subset \complex_g$, and that $\gV^\star$ is {\em outer stabilisable} if a friend $F$ of $\gV^\star$ exists such that $\sigma(A+B\,F\,|\,\gX/\gV^\star) \subset \complex_g$.
 The eigenvalues of $A+B\,F$ restricted to $\gV^\star$ can be further split into two disjoint sets: the eigenvalues of $\sigma(A+B\,F |\gR^\star)$, are all freely assignable\footnote{An assignable set of eigenvalues is always intended to be a set of complex numbers mirrored with respect to the real axis. }
 with a suitable choice of $F$ in $\mathfrak{F}(\gV^\star)$. The eigenvalues of $\sigma\,(A+B\,F | {\gV^\star}/{\gR^\star})$ are fixed for all the choices of $F$ in $\mathfrak{F}(\gV^\star)$ {and coincide with the invariant zeros of $\Sigma$}.
 Thus, $\gV^\star$ is inner stabilisable if and only if $\gZ \subset \complex_g$. 
 Finally, we use the symbol $\gV^\star_g$ to denote the largest stabilisability output-nulling subspace of $\Sigma$, i.e., the largest subspace $\gV$ of $\gX$ for which a matrix $F\,{\in}\,\mathbb{R}^{m\,{\times}\,n}$ exists such that $(A+B\,F)\,\gV\subseteq \gV \subseteq \ker (C+D\,F)$ and $\sigma(A+B\,F\,|\,\gV) \subset \complex_g$. There holds $\gR^\star \subseteq \gV^\star_g \subseteq \gV^\star$. \\[-2mm]

 \subsection{The tracking control problem}

In this paper, we are concerned with the problem of the design of a
state-feedback control law for (\ref{sys}) such that for all initial conditions the output $y$ tracks a step reference $r \in \gY$ with zero steady-state error and is monotonic in all components.
If $y$ asymptotically tracks the constant reference $r$ and is monotonic, then it is also both non-overshooting and non-undershooting. The converse is obviously not true in general.
 The following standing assumption is standard for tracking control problems (see e.g. \cite{He-CW-05}), and ensures that any given constant
reference target $r$ can be tracked from any initial
condition $x_0 \in \gX$:\\[-3mm]
 \begin{assumption}
 \label{Ass1}
System $\Sigma$ is right invertible and stabilisable. {We also assume that
 $\Sigma$ has no invariant
zeros at the origin in the continuous time case, or at $1$ in the discrete case.}\\[-3mm]
\end{assumption}

Right invertibility and the absence of invariant zeros at the origin in the continuous time (or at $1$ in the discrete case) are natural assumptions for a tracking control problem, as they ensure that any given constant reference target $r$ can be tracked from any given initial condition, see \cite{He-CW-05}. These two assumptions generically hold when $m \ge p$.
Under Assumption \ref{Ass1}, the standard method for designing a tracking controller for a
step reference signal is carried out as follows. Given the step reference $r \in \gY$ to track, choose a feedback gain matrix
$F$ such that $A + B\,F$ is asymptotically stable: this is always possible since the pair $(A,B)$ is assumed to be stabilisable. Let us choose two vectors $x_{\rm ss} \in
\gX$ and $u_{\rm ss} \in \gU$ that, for any $r \in \gY$, satisfy 
\begin{eqnarray}
 \label{k2}
\left\{ \begin{array}{rcl}
0 \ns&\ns = \ns&\ns A\,x_{\rm ss}+B\,u_{\rm ss} \\
r \ns&\ns = \ns&\ns C\,x_{\rm ss}+D\,u_{\rm ss} \end{array} \right. \qquad \quad \text{and} \qquad \quad
\left\{ \begin{array}{rcl}
 x_{\rm ss} \ns&\ns = \ns&\ns A\,x_{\rm ss}+B\,u_{\rm ss} \\ 
 r \ns&\ns = \ns&\ns C\,x_{\rm ss}+D\,u_{\rm ss}
\end{array} \right.
\end{eqnarray}
in the continuous and in the discrete case, respectively. 
 Such pair of vectors $x_{\rm ss} \in \gX$ and $u_{\rm ss} \in \gU$ exist since {\em (i)} right invertibility ensures that the system matrix pencil
$P_{\scriptscriptstyle \Sigma}(\lambda)$ is of full row-rank for all but finitely many $\lambda \in \complex$, see \cite[Theorem 8.13]{Trentelman-SH-01}, and, as already recalled, the values $\lambda \in \complex$ for which $P_{\scriptscriptstyle \Sigma}(\lambda)$ loses rank are invariant zeros of $\Sigma$; {\em (ii)} in the continuous (resp. discrete) time case, the absence of invariant zeros at the origin (resp. at $1$) guarantees that the matrix $P_{\scriptscriptstyle \Sigma}(0)$ (resp. $P_{\scriptscriptstyle \Sigma}(1)$) is of full row-rank. As such, Assumption \ref{Ass1} guarantees that the linear system $\left[ \begin{smallmatrix} 0 \\[1mm] r \end{smallmatrix} \right]=\left[ \begin{smallmatrix} A& B \\[1mm] C & D \end{smallmatrix} \right]\left[ \begin{smallmatrix} x_{\rm ss} \\[1mm] u_{\rm ss} \end{smallmatrix} \right]$ in the continuous time or
$\left[ \begin{smallmatrix} 0 \\[1mm] r \end{smallmatrix} \right]=\left[ \begin{smallmatrix} A-I_n& B \\[1mm] C & D \end{smallmatrix} \right]\left[ \begin{smallmatrix} x_{\rm ss} \\[1mm] u_{\rm ss} \end{smallmatrix} \right]$
in the discrete time is always solvable in $\left[ \begin{smallmatrix} x_{\rm ss} \\[1mm] u_{\rm ss} \end{smallmatrix} \right]$. Now, applying the state-feedback control law
 \be
 u(t) = F\,\Big( x(t)-x_{\rm ss} \Big) + u_{\rm ss}
 \label{ulaw}
 \ee
 to (\ref{sys}) and using the change of variable $\xi \defi x -x_{\rm ss}$ gives the closed-loop
 homogeneous system
 \be
 \Sigma_{hom}: \
\left\{ \begin{array}{lcr}
\gD\,{\xi}(t) \ns&\ns = \ns&\ns (A+B\,F)\,\xi(t),\;\;\quad \xi(0)=x_0-x_{\rm ss},\hfill\cr
 y(t) \ns&\ns = \ns&\ns (C+D\,F)\,\xi(t)+ r. \hfill \end{array} \right.
 \label{syschom}
 \ee
Since $A+B\,F$ is asymptotically stable, $x$ converges to $x_{\rm
ss}$, $\xi$ converges to zero and $y$ converges to $r$ as $t$ goes to infinity. We shall refer to $\xi$ as the {\em error state coordinates}.

 \subsection{Achieving a globally monotonic response with any desired convergence rate}

In this paper we are concerned with the general problem of finding a gain matrix $F$ such that the closed-loop system obtained using \eqref{ulaw} in \eqref{sys} achieves a monotonic response at any desired rate of convergence, from all initial conditions. We shall describe this property as {\em global monotonicity}. We describe the problem as follows in terms of the tracking error $\epsilon(t)\defi y(t)-r(t) = \bsmat \epsilon_1(t) \\[-1mm] \vdots \\[1mm] \epsilon_p(t)\esmat$.

\begin{problem} (Global Monotonicity) 
\label{prob:mono0}
Let $\rho \in \real$, such that $\rho<0$ in the continuous time and $\rho \in (0,1)$ in the discrete time. Find a state-feedback matrix $F$ such that applying (\ref{ulaw}) with this $F$ to $\Sigma$ yields an asymptotically stable closed-loop system $\Sigma_{hom}$ for which the tracking error term $\epsilon(t)$ converges monotonically to 0 at a rate at least $\rho$ in all outputs, from all initial conditions. Specifically, we require that in the continuous time
\begin{equation}
 \forall \xi(0) \in \gX,\	\forall k \in \{1, \dots, p\}, \ \exists \,\beta_k \in \mathbb{R} \, : \;\; | \epsilon_k(t) | \leq \beta_k \, \exp (\rho\, t) \quad \forall t \in \real_+
\label{eq:error2}
\end{equation}
where $\epsilon_k(t)$ is strictly monotonic in $t$, and in the discrete time
\begin{equation}
 \forall \xi(0) \in \gX,\	\forall k \in \{1, \dots, p\}, \ \exists \,\beta_k \in \mathbb{R} \, : \;\; | \epsilon_k(t) | \leq \beta_k \, {\rho^{t}}\quad \forall t \in \mathbb{N},
\label{eq:error2d}
\end{equation}
where, again, $\epsilon_k(t)$ is strictly monotonic in $t$.
\end{problem}


If we are able to obtain a tracking error $\epsilon(t)$ that consists of a single exponential per component in the continuous time or a single power per component in the discrete time, i.e., 
\begin{equation} 
\label{eq:error}
\epsilon(t) =
 \begin{bmatrix} \beta_1\,\exp(\lambda_1\,t)\\[-2mm]\vdots \\[-1mm] \beta_p\,\exp(\lambda_p\,t) \end{bmatrix} \qquad \text{or} \qquad
\epsilon(t) =
 \begin{bmatrix} \beta_1\,\lambda_1^{t} \\[-2mm]\vdots \\[-1mm] \beta_p\,\,\lambda_p^{t} \end{bmatrix},
 \end{equation}
 respectively, and we are able to choose each $\lambda_k$ in such a way that $\lambda_k \leq \rho<0$ in the continuous time and $0 \le \lambda_k \leq \rho$ in the discrete time, 
then we obtain a solution to Problem \ref{prob:mono0}. Indeed, asymptotically stable exponentials of $\lambda_k$ or powers of $\lambda_k$ are monotonic functions. This suggests that a possible way of solving Problem \ref{prob:mono0} consists in the solution of the following problem.


\begin{problem} (Single mode outputs).
\label{prob:mono1}
Let $\rho \in \real$, such that $\rho<0$ in the continuous time and $\rho \in (0,1)$ in the discrete time. Find a feedback matrix $F$ and a $p$-tuple of distinct values $\mathfrak{L}_p =(\lambda_1, \lambda_2, \cdots, \lambda_p)$ 
such that, in the continuous time $\lambda_k \leq \rho<0$, and in the discrete time $0 \le \lambda_k \leq \rho$ 
such that applying (\ref{ulaw}) with this $F$ to $\Sigma$ yields an asymptotically stable closed-loop system $\Sigma_{hom}$ for which, from all initial conditions, the tracking error term is given by
(\ref{eq:error})
for some real coefficients $\{\beta_k \}_{k=1}^p$ depending only upon $\xi(0)$.
\end{problem}


If we can guarantee that $F$ yields (\ref{eq:error}) for any initial condition and any $r \in \gY$, and such that $\lambda_k \leq \rho<0$ in the continuous time and $0 \le \lambda_k \leq \rho$ in the discrete time, then obviously such feedback also solves Problem \ref{prob:mono0}. However, the following result shows that the converse is true as well, i.e., the only way to obtain a feedback that ensures global monotonic tracking with arbitrary rate of convergence of the tracking error is to obtain a tracking error in the form (\ref{eq:error}).

\begin{lemma} 
\label{lem:mono}
Problem~\ref{prob:mono0} is equivalent to Problem~\ref{prob:mono1}.
\end{lemma}
\proof
Let us consider the continuous time case, the discrete case being entirely equivalent.
 Let $\rho$ be an arbitrary negative real number. It is clear that if $F$ and $\mathfrak{L}_p$ solve Problem~\ref{prob:mono1} with respect to this $\rho$, then the outputs $\epsilon_k(t)$ satisfy (\ref{eq:error}), and hence also (\ref{eq:error2}) or (\ref{eq:error2d}). 
 Next, assume that the feedback matrix $F$ solves Problem~\ref{prob:mono0} for a certain $\rho \in \real^-$. For some $1 \leq \nu \leq n$, let $\mathcal{L} = \{\lambda_1,\lambda_2, \ldots, \lambda_\nu \}$ be the eigenvalues of $A + B\,F$, with associated algebraic multiplicities $\mathcal{M} = \{m_1,m_2, \ldots, m_\nu\}$ satisfying $m_1 + m_2+ \ldots + m_\nu= n$.
If $\mathcal L$ contains any complex eigenvalues we shall assume these are ordered such that $\lambda_i = \overline{\lambda}_{i+1}$. Applying (\ref{ulaw}) with this $F$ to $\Sigma$, we obtain $\Sigma_{hom}$ in
(\ref{syschom}). Let $V \defi [\,V_1 \;\; V_2 \;\;\dots\;\; V_\nu\, ]$ denote the eigenvector matrix of $A+B\,F$, where $V_i$ denotes a column matrix of $m_i$ generalised eigenvectors associated with $\lambda_i$. Then $V$ is non-singular, and for any initial condition $\xi_0 \defi \xi(0) \in \gX$ of (\ref{syschom}), we can introduce
 $\alpha = [ \alpha_1 \ \alpha_2 \ \ldots \ \alpha_\nu]^\top = V^{-1} \xi(0)$, where each $\alpha_i$ is a column matrix of dimension $m_i$ associated with $V_i$. For the case of $\lambda_i$ complex, we have $V_i = \overline{V}_{i+1}$ and hence also $\alpha_{i} = \overline{\alpha}_{i+1}$.
 The tracking error arising from $\xi_0$ is given by 
 \be
 \label{error}
\epsilon(t) = -(C+D\,F)\,\exp[(A+B\,F)\,t]\,\xi_0 = (C+D\,F) \, \sum_{i=1}^{\nu} \sum_{j=1}^{m_i} \gamma_{i,j}\, t^{j-1} \exp(\lambda_i t),
\ee
where $\gamma_{i,j}$ is a function of the entries of $V_i$ and $\alpha_i$. Thus, each component of the tracking error is comprised of functions of the form
\be 
t^{j-1}\,\exp(\lambda_i\, t), \;\; t^{j-1} \exp(\sigma_i\, t)\cos(\omega_i\, t), \;\; t^{j-1} \exp(\sigma_i \,t)\sin(\omega_i\, t) 
\ee
with $0 \le j-1 \le m_i$, and 
where for complex $\lambda_i$ we have $\sigma_i = \mathfrak{Re}\{\lambda_i\}$ and $\omega_i = \mathfrak{Im}\{\lambda_i\}$. By (\ref{eq:error2}), we conclude that $\mathfrak{Re}\{\lambda_i\} \leq \rho$ for all $\lambda_i \in \mathcal L$. Further, by assumption the response is monotonic, which implies that only components of the form $\exp(\lambda_i \,t)$ for real $\lambda_i$ can appear in each output (notice that, for any real $\lambda_i<0$, the function $t^{j-1}\,\exp(\lambda_i\, t)$ is monotonic only if $j =1$). Thus for each output $\epsilon_k(t)$ we have real coefficients $\tilde{\beta}_{k,i}$, depending on $\xi_0$, such that
\be
\epsilon_k(t) = \sum_{i=1}^{\nu} \tilde{\beta}_{k,i} \exp(\lambda_i t).
\ee
Moreover each response is monotonic from all initial conditions. From Lemma A.1 of \cite{Schmid-N-09}, if $\epsilon_k(t)$ is the sum of two or more negative real exponential functions, it will change sign (and hence not be monotonic) for some values of the coefficients $\tilde{\beta}_{k,i}$. Since the mapping $\xi_0 \mapsto \alpha $ is surjective, the response can only be globally monotonic if each output component is comprised of a single real exponential function. Thus, we must have
$\epsilon_k(t) = \tilde{\beta}_{k,i} \exp(\lambda_i t)$
for some $\lambda_i \in \gL$ and some real coefficient $\beta_k \defi \tilde{\beta}_{k,i}$. If we select the $k$-th element of the $p$-tuple $\mathfrak{L}_p$ to equal the element of the set $\mathcal L$ that
appears in the tracking error component $\epsilon_k$, Problem 2 is solved. 
\endproof
\ \\[-4mm]
Another important and useful problem is one in which the requirements include a specified choice of the closed-loop modes that are visible in each component of the tracking error:
\ \\[-4mm]
\begin{problem} (Single mode outputs).
\label{prob:mono2}
Let $\mathfrak{L}_p =(\lambda_1, \lambda_2, \cdots, \lambda_p)$ be a $p$-tuple of real numbers such that, in the continuous time $\lambda_i <0$ and in the discrete time $\lambda_i \in (0,1)$ for all $ i \in \{1,\ldots,p\}$. Find a feedback matrix $F$ such that applying (\ref{ulaw}) 
to $\Sigma$ yields an asymptotically stable closed-loop system $\Sigma_{hom}$ for which, from all initial conditions and for all step references, the tracking error term is given by (\ref{eq:error}). 
\end{problem}

\section{Global monotonicity: the intuitive idea} 
\label{sec:intu}
In the previous section, we observed that in order for the problem of global monotonic tracking to be solvable, we need to render at least $n-p$ closed-loop modes invisible from the tracking error and distribute the remaining $p$ modes evenly into the tracking error with one mode per error component. If this is possible, then the step response is guaranteed to be monotonic for any initial condition, and therefore also non-overshooting and non-undershooting. 
The converse is true as well, as shown in Lemma \ref{lem:mono}.
If we are able to render more than $n-p$ modes invisible at $\epsilon(t)$, one or more components of the tracking error can be rendered identically zero, and therefore for those components instantaneous tracking can also be achieved, in which the output component immediately takes the desired reference value.
The aim of the next part of the paper is to find conditions under which a gain matrix $F$ can be obtained to deliver the single mode structure (\ref{eq:error}) for any initial condition. 
Consider $\Sigma_{hom}$ in (\ref{syschom}), which can be re-written as
 \be
 \Sigma_{hom}: \
\left\{ \begin{array}{lcr}
\gD\,{\xi}(t) \ns&\ns = \ns&\ns A\,\xi(t)+B\,\omega(t),\;\;\; \hfill\cr
\epsilon(t) \ns&\ns = \ns&\ns C\,\xi(t)+ D\,\omega(t), \hfill \end{array} \right.
 \label{syschom1}
 \ee
where $\omega(t)=F\,\xi(t)$. Clearly, $\Sigma_{hom}$ can be identified with the quadruple $(A,B,C,D)$. 
The task is now to 
find a feedback matrix $F$ such that the new control $\omega(t)=F\,\xi(t)$ guarantees that for every initial condition $\xi_0\in\gX$ the tracking error $\epsilon(t)$ is characterised by a single stable real mode per component.
Let $j \in \{1,\ldots,p\}$. Let $\lambda_j$ be real, stable and not coincident with any of the invariant zeros of $\Sigma$. Consider a solution $v_j$ and $w_j$ of the linear equation
\bea
\label{ML0}
\left[ \begin{array}{ccc} \! A-\lambda_j\,I_n \! & \! B \! \\ \! C \! & \! D \! \end{array} \right]\left[ \begin{array}{ccc} v_j \\ w_j \end{array} \right]=\left[ \begin{array}{ccc} 0 \\ \beta_j\,e_j \end{array} \right],
\eea
where $\beta_j \neq 0$ and $e_j$ is the $j$-th vector of the canonical basis of $\gY$. Notice that (\ref{ML0}) always has a solution in view of the right-invertibility of $\Sigma$. By choosing $F$ such that $F\,v_j=w_j$, we find $(A+B\,F)\,v_j=\lambda_j\,v_j$ and $(C+D\,F)\,v_j=\beta_j\,e_j$. Hence, from (\ref{error}) we know that for any initial error state $\xi_0\in\spanR\{v_j\}$ the response associated with the control $\omega(t)=F\,\xi(t)$ is
\bea
\label{newalpha}
\epsilon(t)=\left[ \begin{array}{ccc} 0 \\[-2mm] \vdots \\[-2mm] \gamma_j\,\exp(\lambda_j\,t) \\[-3mm] \vdots \\[-2mm] 0 \end{array} \right] \left. \begin{array}{ccc} \phantom{0} \\[-2mm] \phantom{\vdots} \\[-2mm] \leftarrow j \phantom{q\,\exp(\lambda\,t)} \\[-3mm] \phantom{\vdots} \\[-2mm] \phantom{0} \end{array} \right.
\eea
where $\gamma_j$ depends on the particular initial state $\xi_0$. Considering $\mathfrak{L}_p=(\l_1,\ldots,\l_p)$ with each $\l_j$ real and stable and different from the invariant zeros, 
by applying this argument for all components of the tracking error, we obtain a set of solutions $\left[\begin{smallmatrix} v_1 \\[1mm] w_1 \end{smallmatrix} \right], \left[\begin{smallmatrix} v_2 \\[1mm] w_2 \end{smallmatrix} \right], \ldots,\left[\begin{smallmatrix} v_p \\[1mm] w_p \end{smallmatrix} \right]$ of (\ref{ML0}). 
If $v_1,\ldots,v_p$ are linearly independent, we can choose
 $F$ to be such that $F\,v_i=w_i$ for all $i \in \{1,\ldots,p\}$. Then, for every $\xi_0 \in \spanR \{v_1,v_2,\ldots,v_p\}$, by superposition we find
\bea
\label{risp}
\epsilon(t)=\left[ \begin{array}{ccc} \gamma_1\,\exp(\lambda_1\,t) \\[-2mm] 0 \\[-3mm]\vdots \\[-2mm] 0 \end{array} \right] +
\left[ \begin{array}{ccc} 0 \\[-2mm] \gamma_2\,\exp(\lambda_2\,t) \\[-3mm]\vdots \\[-2mm] 0\end{array} \right] +\ldots+
\left[ \begin{array}{ccc} 0 \\[-2mm] 0 \\[-2mm]\vdots \\[-2mm] \gamma_p\,\exp(\lambda_p\,t) \end{array} \right] =
\left[ \begin{array}{ccc} \gamma_1\,\exp(\lambda_1\,t) \\[-1mm] \gamma_2\,\exp(\lambda_2\,t) \\[-2mm]\vdots \\[-1mm] \gamma_p\,\exp(\lambda_p\,t) \end{array} \right].
\eea
However, this result only holds when $\xi_0 \in \spanR \{v_1,v_2,\ldots,v_p\}$.
In order for this response to be achievable from any initial condition, we also need to render  the remaining $n-p$ closed-loop modes invisible at $\epsilon(t)$. This task can be accomplished by exploiting the supremal stabilisability output-nulling subspace $\gV^\star_g$ of the system, which is defined as the largest subspace of $\gX$ for which a friend $F$ exists such that, for every initial state lying on it, the corresponding state feedback generates a state trajectory that asymptotically converges to zero while the corresponding output (the tracking error in the present case) remains  at zero.
We will see in Section \ref{background} that a basis for $\gV^\star_g$ can always be obtained 
as the image of a matrix $[\,V_1\;\;\;V_2\;\;\ldots \;\; V_d\,]$ that satisfies
\bea
\label{RosVg}
\left[ \begin{array}{ccc} \! A-\mu_j\,I_n \! & \! B \! \\ \! C \! & \! D \! \end{array} \right]\left[ \begin{array}{ccc} V_{j} \\ W_{j} \end{array} \right]=0,
\eea
for some other matrix $[\,W_1\;\;\;W_2\;\;\ldots \;\; W_d\,]$ partitioned comformably, 
where $\{\mu_1,\ldots,\mu_t\}$ are the (distinct) minimum-phase invariant zeros of $\Sigma$ and $\{\mu_{t+1},\ldots,\mu_d\}$ are arbitrary and stable (let us assume for the moment that they are real and distinct).
If the dimension of $\gV^\star_g+\spanR \{v_1,v_2,\ldots,v_p\}$ is equal to $n$, every initial state $\xi_0\in\gX$ can be decomposed as the sum $\xi_v+\xi_r$, where $\xi_v\in \gV^\star_g$ and $\xi_r \in \spanR \{v_1,v_2,\ldots,v_p\}$. 
 If for the sake of argument we have $\dim \gV^\star_g=n-p$, 
 and we can find a set of linearly independent column vectors $\{v_{p+1},\ldots,v_n\}$ from the columns of $[\,V_1\;\;\;V_2\;\;\ldots \;\; V_d\,]$ that is linearly independent of $\{v_1,\ldots,v_p\}$, we can take $w_{p+1},\ldots,w_{n}$ to be the columns of $W_g$ that correspond to $v_{p+1},\ldots,v_{n}$, and construct the feedback control $\omega(t)=F\,\xi(t)$ where $F$ is such that
 $F\,[\,v_1\;\; \ldots \;\; v_p \;\; v_{p+1} \;\; \ldots \;\; v_{n}\,]=[\,w_1\;\; \ldots \;\; w_p \;\; w_{p+1} \;\; \ldots \;\; w_{n}\,]$, the response associated with $\xi_v$ is identically zero, while the one associated with $\xi_r$ is still given by (\ref{risp}). Hence, the tracking error can be written as in (\ref{eq:error}) for any $\xi_0 \in \gX$. The closed-loop eigenvalues obtained with this matrix $F$ are given by the union of $\{\l_1,\ldots,\l_p\}$, with the set of values $\mu_j$ that are associated with the columns $\{v_{p+1},\ldots,v_n\}$ chosen from $[\,V_1\;\;\;V_2\;\;\ldots \;\; V_d\,]$.


We now provide some more intuition on our design method by using an example that embodies all those system theoretic characteristics that are perceived as the major difficulties in achieving monotonic tracking. This system is MIMO, bi-proper, uncontrollable (although obviously stabilisable) and is characterised by 3 non-minimum phase zeros. To the best of the authors' knowledge, there are no methods available in the literature that can solve the tracking problem for MIMO systems with a guaranteed monotonic response under such assumptions, and especially in the presence of three non-minimum phase invariant zeros. 
We also want to stress that this problem is even solved in closed form.

\begin{example}
\label{exe0}
{Consider the bi-proper continuous-time LTI system $\Sigma$ in (\ref{sys}) with
\beann
A \ns&\ns = \ns&\ns \bsmat
 -6 && 0 && 0 && 0 && 0\\[1mm]
 3 && 3 && 0 && 0 && 0\\[1mm]
 0 && 0 && 2 && 0 && 2\\[1mm]
 -1 && 0 && 2 && 0 && 0\\[1mm]
 -2 && 0 && 0 && 0 && 2\esmat\!, \quad B=\bsmat
 0 && 0 && 0 && 0\\[1mm]
 0 && 0 && 0 && -3\\[1mm]
 0 && 4 && 2 && 0\\[1mm]
 1 && -1 && 0 && -1\\[1mm]
 0 && -1 && 0 && 0 \esmat\!,\quad 
 C = \bsmat
 -1 && 0 && 0 && 0 && 0\\[1mm]
 3 && 0 && 0 && 0 && 9\\[1mm]
 1 && 0 && 0 && 0 && 0\esmat\!, \qquad D=\bsmat
 0 && 0 && -2 && 0\\[1mm]
 0 && 3 && -3 && -3\\[1mm]
 0 && 0 && 2 && -2 \esmat.
 \eeann
 We want to find a feedback matrix $F$ such that the output of this system monotonically tracks a unit step in all output components, and the assignable closed-loop eigenvalues are equal to $\lambda_1=-1$, $\lambda_2=-2$ and $\lambda_3=-1$ for the corresponding error components.
 This system is not square since the number of inputs exceeds the number of outputs. It is seen to be right invertible but not left invertible. 
 This system is not reachable but it is stabilisable, since the only uncontrollable eigenvalue is $-6$, and is equal to the only minimum-phase invariant zero $\mu_1=z_1=-6$ of the system. 
 The other three invariant zeros of the system are non-minimum-phase, and their values are $z_2=2$, $z_3=3$ and $z_4=5$. 
 The first $p=3$ modes must be evenly distributed among the $3$ components of the tracking error. The subspace of initial conditions for which the closed-loop mode governed by the eigenvalue $\lambda_j$ appears in the $j$-th output-component is given by the span of $v_j$, where $v_j$ solves (\ref{ML0})
 for a suitable $w_j\in \gU$. The linear equation (\ref{ML0}) can be solved by pseudo-inversion with $\lambda_1=-1$, $\lambda_2=-2$ and $\lambda_3=-1$, and gives 
 \beann
 \bsmat v_1 \\[1mm] w_1 \esmat \ns&\ns = \ns&\ns P_{\scriptscriptstyle \Sigma}^\dagger(\lambda_1) \bsmat 0 \\[1mm] e_1 \esmat= \frac{1}{18}
 [\, 0 \;\;\; -{27}/{4} \;\;\; 20 \;\;\; -29 \;\;\; -3 \;\;\; | \;\;\;
 -29 \;\;\; -9 \;\;\; -9 \;\;\; -9\,]^\top,\\
 \bsmat v_2 \\[1mm] w_2 \esmat \ns&\ns= \ns&\ns P_{\scriptscriptstyle \Sigma}^\dagger(\lambda_2) \bsmat 0 \\[1mm] e_2 \esmat= \frac{1}{21}[\, 0 \;\;\; 0 \;\;\; -{9}/{2} \;\;\; {26}/{5} \;\;\; 1\;\;\; | \;\;\;
 {13}/{5} \;\;\; 4 \;\;\; 0 \;\;\; 0\,]^\top, \\
 \bsmat v_3 \\[1mm] w_3 \esmat \ns&\ns= \ns&\ns P_{\scriptscriptstyle \Sigma}^\dagger(\lambda_3) \bsmat 0 \\[1mm] e_3 \esmat= \frac{1}{18} [\, 0 \;\;\; -{27}/{4} \;\;\; 7 \;\;\; -{55}/{4} \;\;\; -{3}/{2}\;\;\; | \;\;\;
 -{55}/{4} \;\;\; -{9}/{2} \;\;\; 0 \;\;\; -9 \,]^\top.
 \eeann

 The subspace $\gV^\star_g$ is spanned by a matrix $V_g$ obtained from a basis matrix $\bsmat V_g \\[1mm] W_g \esmat$ for the null-space of $\ker \bsmat A-z_1\,I_n & B \\[1mm] C & D \esmat$. In this specific case we find
 $V_g=\bsmat -2 && 2/3 && -41/22 && 0 && -1/11 \\[1mm] 0 && 0 && 0 && 1 && 0 \esmat^\top$ and 
 $W_g=\bsmat 5 && 36/11&& 1 && 0 \\[1mm] -6 && 0 && 0 && 0 \esmat^\top$.
 Indeed, the subspace of $\gX$ spanned by any basis matrix of the null-space of $P_{\scriptscriptstyle \Sigma}(\lambda)$ for any other value $\lambda$ is linearly dependent of $\ima V_g$, and therefore does not contribute to the construction of a basis for $\gV^\star_g$. 
 Since the dimension of $\gV^\star_g$ is equal to $2$, two closed-loop modes can be rendered invisible at the tracking error. 
 The columns of the square and non-singular matrix $V=[\, v_1 \;\; v_2 \;\; v_3 \;\;\; V_g\,]$ span the subspace $\gV^\star_g+\spanR\{v_1,v_2,v_3\}$, which is $n$-dimensional. 
 The response from any arbitrarily assigned initial error state $\xi_0 \in \gX$ is given by the sum of the responses obtained by projecting the initial state on the subspaces $\ima V_g=\gV^\star_g$, $\spanR \{v_1\}$, $\spanR \{v_2\}$ and $\spanR \{v_3\}$. The first projection gives a null contribution, because the corresponding state trajectory is invisible at the tracking error by definition of $\gV^\star_g$. The projection of the initial state on $\spanR \{v_j\}$ is such that only the $j$-th component of the response is affected and is governed by the mode $\exp(\lambda_j\,t)$ as in (\ref{newalpha}), while all other components are equal to zero. By superposition, the overall state trajectory originating from $\xi_0$ yields a tracking error with one mode per component. The feedback matrix is given in closed form by
 \bea
 \label{f1}
 F=W\,V^{-1}=\left[ \begin{array}{ccccccc}
 \frac{68419}{8250} & \frac{802}{125} & -\frac{1121}{125} & -6 & -\frac{1639}{250}\\
 -\frac{5351}{2475} & -\frac{16}{75} & \frac{6}{25} & 0 & \frac{127}{25}\\
 \frac{5537}{4950} & -\frac{12}{225} & -\frac{36}{25} & 0& -\frac{162}{25}\\
 \frac{4}{9}& \frac{4}{3} & 0 & 0 & 0 \end{array} \right],
 \eea
 where $W=[\,w_1 \;\; w_2 \;\; w_3 \;\;\, W_g\,]$. The closed-loop eigenvalues are $\sigma(A+B\,F)=\{-2,-1,-6\}$ where the multiplicity of the eigenvalues $-1$ and $-6$ is equal to two as expected. In fact, these closed-loop eigenvalues are given by the union of the minimum-phase invariant zeros and the assigned closed-loop eigenvalues. The closed-loop eigenvalue corresponding to the minimum-phase zero appears with double multiplicity because the kernel of the Rosenbrock matrix $P_{\scriptscriptstyle \Sigma}(\lambda)$ evaluated at that minimum-phase invariant zero spans a two dimensional subspace. 
 If the reference is $r=[2\;\;\; \,2\;\;\;\, 2]^\top$, we compute $x_{\rm ss}$ and $u_{\rm ss}$ by solving (\ref{k2}), and we obtain $x_{\rm ss} =[ 0\;\; -2 \;\; {10}/{3} \;\; 0 \;\; -{7}/{15}\,]^\top$ and $u_{\rm ss}=[\, -{48}/{5} \;\;
 -{14}/{15}\;\; -1\;\; -2 \,]^\top$.
 Given an arbitrary initial condition $\xi_0 \in \gX$, we compute $\alpha=[\,\alpha_1\;\; \ldots\;\;\alpha_5\,]^\top$ from $\alpha=V^{-1}\,\xi_0$ and the tracking error that follows from the application of the control law $u(t) = F\,( x(t)-x_{\rm ss} ) + u_{\rm ss}$ with the feedback matrix $F$ in (\ref{f1}), yields the tracking error
 $\epsilon(t)=[\, \gamma_1\,\exp(-t) \;\; \gamma_2\,\exp(-2\,t) \;\; \gamma_3\,\exp(-t) \,]^\top$, which has the single mode form of (\ref{risp}). Therefore, the system exhibits a globally monotonic step response.
 The tracking errors of the closed-loop system are shown in Figure \ref{fig1} for two different initial conditions. 
\begin{figure}
\hspace{1cm}\includegraphics[width=6.5cm, height=3.5cm]{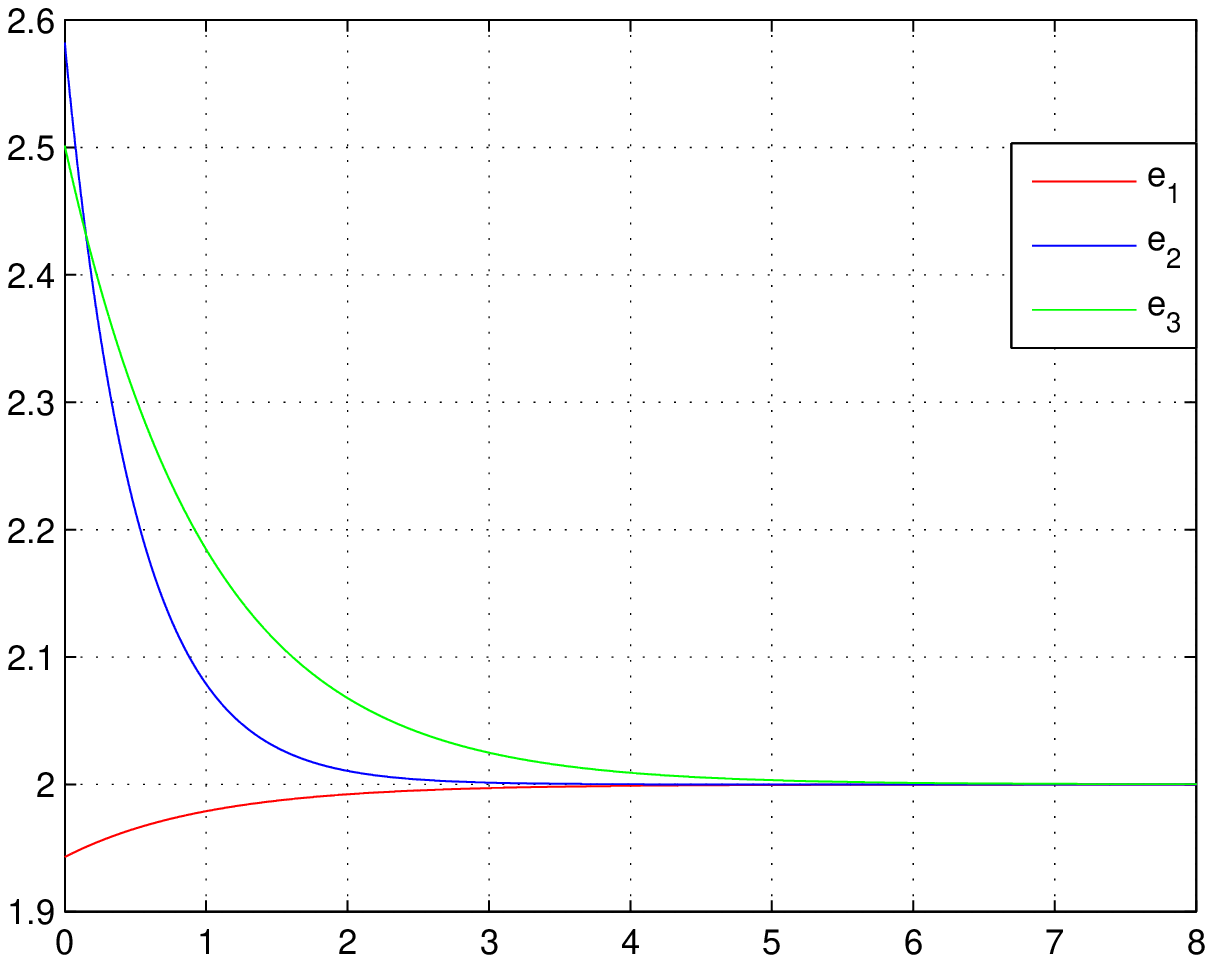}
\includegraphics[width=6.5cm, height=3.5cm]{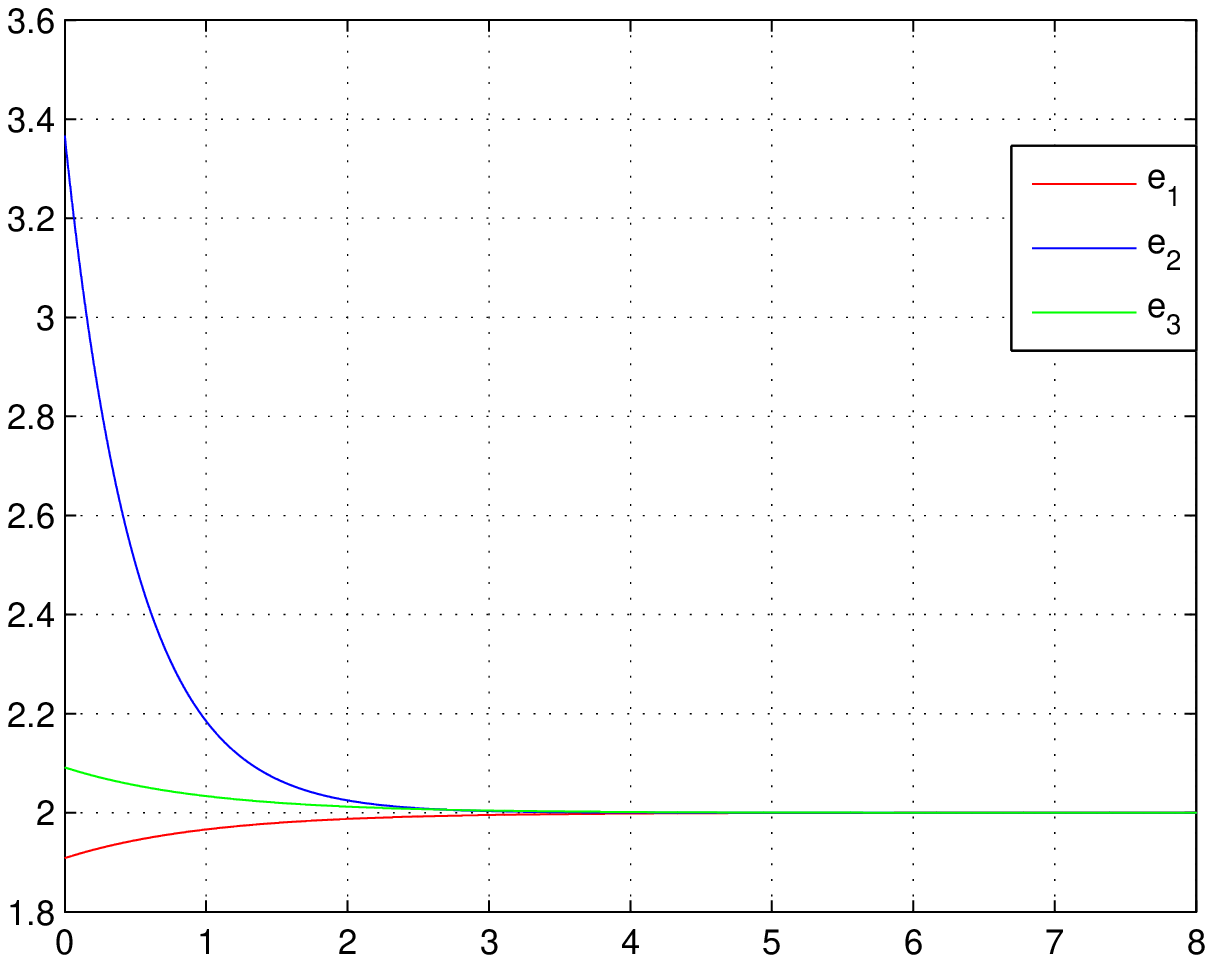}
\caption{Tracking errors of $\Sigma$ with respect to the reference $r=[\,2\;\;\;2\;\;\; 2\,]^\top$ with initial conditions $x_0=[\,0.1\;\; -\!0.2\;\, 0.1\;\, 0.1\;\, 0\,]^\top$ and $x_0=[\,0.6 \,\; 0.2\,\; 0.2\;\, -\!0.2\,\; 1\,]^\top$, respectively.}
 \label{fig1}
 \end{figure}
 \endproof
 }
 \end{example}

Since the solvability condition for global monotonicity is given in terms of the dimension of the subspace $\gV^\star_g+\spanR \{v_1,\ldots,v_p\}$, and that $v_1,\ldots,v_p$ depend on the choice of the closed-loop eigenvalues $\lambda_1,\ldots,\lambda_p$, the solvability condition given in the previous example seems to depend on the particular choice of the closed-loop eigenvalues. 
The question at this point is: how does the choice of the closed-loop eigenvalues affect the dimension of $\gV^\star_g+\spanR \{v_1,\ldots,v_p\}$? Are there good and bad choices of the closed-loop eigenvalues? More generally, can we find alternative solvability conditions given solely in terms of the system structure and not in terms of a choice of eigenvalues? These are the crucial points that will be addressed in the sequel.

\section{Mathematical background}
\label{background}
%
In the previous section we identified the basic tools 
that can be used to obtain a gain matrix $F$ such that $p$ of the $n$ closed-loop modes are evenly distributed into the $p$ components of the tracking error as in (\ref{risp}), and the remaining $n-p$ modes are rendered invisible at the tracking error.
%
The first of these tools is the subspace $\gV^\star_g$,
which is made up of the sum of two parts. The first is the subspace $\gR^\star$, and the second is, loosely, the subspace spanned by the directions of the minimum-phase invariant zeros of $\Sigma$. In this section, we recall some important results concerning the relations between these subspaces and the null-space of the Rosenbrock system matrix pencil $P_{\scriptscriptstyle \Sigma}$.
 %
%
Given $\mu \in \complex$, we use the symbol $N_{\scriptscriptstyle \Sigma}(\mu)$ to denote a basis matrix for the null-space of $P_{\scriptscriptstyle \Sigma}(\mu)$, and we denote by $d(\mu)$ the dimension of this null-space. Let $d \defi n+m-\normrank P_{\scriptscriptstyle \Sigma}(\mu)$. There holds $d(\mu)=d$, unless {$\mu\in\gZ$, i.e., $\mu$} is an invariant zero of $\Sigma$, in which case $d(\mu)>d$.
Given a set of $h$ self-conjugate complex numbers $\gL = \{\mu_1,\ldots,\mu_h\}$ containing exactly $s$ complex conjugate pairs, we say that $\gL$ is $s$-conformably ordered if $2\,s \le h$ and the first $2\,s$ values of $\gL$ are complex, while the remaining are real, and for all odd $k \le 2\,s$ we have $\mu_{k+1}=\overline{\mu}_k$. For example, the sets $\gL_1=\{1+\i,1-\i,3,-4\}$, $\gL_2=\{10\,\i,-10\,\i,2+2\,\i,2-2\,\i,7\}$ and $\gL_3=\{3,-1\}$ are respectively $1$-, $2$- and $0$-conformably ordered. \\[0mm]

\subsection{Computation of a basis of $\gR^\star$}

The following result, see \cite{SICON} and \cite{CDC}, presents a procedure for the computation of a basis matrix for $\gR^\star$ and, simultaneously, for the parameterisation of all the friends of $\gR^\star$ that place the eigenvalues of the closed-loop restricted to $\gR^\star$ at arbitrary locations. This procedure aims at constructing a basis for $\gR^\star$ starting from basis matrices $N _{\scriptscriptstyle \Sigma}(\mu_i)$ of the null-spaces of the Rosenbrock matrix relative to an $\sigma$-conformably ordered set $\gL= \{\mu_1,\ldots,\mu_r\}$, where $r\defi \dim \gR^\star$, which will result as closed-loop eigenvalues. No generality is lost by assuming that for every odd $i \in \{1,\ldots,2\,\sigma\}$, the basis matrix $N _{\scriptscriptstyle \Sigma}(\mu_{i+1})$ is constructed as $N _{\scriptscriptstyle \Sigma}(\mu_{i+1})=N _{\scriptscriptstyle \Sigma}(\overline{\mu_{i}})=\overline{N _{\scriptscriptstyle \Sigma}(\mu_{i})}$.\\[0mm]


\begin{lemma}{\sc (\cite{SICON,CDC})}.
\label{Rs}
Let $r = \dim \gR^\star$. Let $\gL = \{\mu_1,\ldots,\mu_r\} \subset \real \setminus \gZ$ be distinct. Let $k_i \in \real^d$ for each $i \in \{1,\ldots,r\}$, and define $K \defi \diag \{k_1,\ldots,k_r\}$. 
Let 
\bea
\label{MK}
\bmat{c} V_{\scriptscriptstyle K}\\[-1mm] W_{\scriptscriptstyle K} \emat \defi \bmat{c|c|c|c} N _{\scriptscriptstyle \Sigma}(\mu_1) & N _{\scriptscriptstyle \Sigma}(\mu_2) & \ldots & N _{\scriptscriptstyle \Sigma}(\mu_r) \emat\,K,
\eea
where $V_{\scriptscriptstyle K} \in \real^{n \times r}$ and $W_{\scriptscriptstyle K} \in \real^{m \times r}$. Then, {\bf (i)} Matrix $V_{\scriptscriptstyle K}$ is generically full column-rank with respect to $K$, i.e., $\rank\, V_{\scriptscriptstyle K} =r$ for every $K$ except for those lying in a set of Lebesgue measure zero; {\bf (ii)} For all $K$ such that $\rank\, V_{\scriptscriptstyle K} =r$, we have $\gR^\star=\ima V_{\scriptscriptstyle K}$; {\bf (iii)} The set of all friends of $\gR^\star$ such that $\sigma(A+B\,F\,|\,\gR^\star)=\gL$ is parameterised in $K$ as $F=W_{\scriptscriptstyle K}\,V_{\scriptscriptstyle K}^{\dagger}$, where $K$ is such that $\rank\,V_{\scriptscriptstyle K} =r$. \\[-4mm]
\end{lemma}

Lemma \ref{Rs} permits us to write a spanning set of $\gR^\star$ in terms of the selection of at most $r$ real numbers, as we now show.
For any $\mu \in \real \setminus \gZ$, let us define
\begin{equation} \label{eq:RsjL}
	\gR^\star (\mu) \stackrel{\text{\tiny def}}{=} \left\{ v \in \gX \,\Big|\,\, \exists w \in \gU \;:\;\;
\left[ \begin{smallmatrix} A-\mu\,I_n & B \\[1mm] C & D \end{smallmatrix} \right]\left[ \begin{smallmatrix} v \\[1mm] w \end{smallmatrix} \right]=0 \right\}.
\end{equation}
From this definition, a decomposition of $\gR^\star$ can be obtained from Lemma~\ref{Rs}. 

\begin{corollary}
\label{agg}
Given any distinct set $\{\mu_{1},\cdots,\mu_{r}\} \subset \real \setminus \gZ$, there
holds 
\bea
\label{decomR}
\gR^{\star}=\gR^{\star}(\mu_{1})+\cdots+\gR^{\star}(\mu_{r}).
\eea
\end{corollary}
\proof 
Let $\bsmat V_i \\[1mm] W_i \esmat$ be a basis matrix of $\ker P_{\scriptscriptstyle \Sigma}(\mu_i)$, and let it be partitioned conformably with $P_{\scriptscriptstyle \Sigma}(\mu_i)$. Matrix $V_i$ is of full column-rank. 
Indeed, if $\omega \in \ker V_i$, then $P_{\scriptscriptstyle \Sigma}(\mu_i)\bsmat V_i \\[1mm] W_i \esmat\,\omega=0$ implies $\bsmat B \\[1mm] D \esmat W_i\,\omega=0$. Since $\bsmat B \\[1mm] D \esmat$ is assumed to be of full column-rank, we conclude that $\omega \in \ker W_i$, so that $\omega\in \ker \bsmat V_i \\[1mm] W_i \esmat=\{0\}$, thus $\omega$ is zero. The matrix $V_i$ is a basis matrix for $\gR^{\star}(\mu_{i})$.
\endproof

\subsection{Computation of a basis of $\gV^\star_g$}
We now turn our attention to the computation of $\gV^\star_g$. From now on, we will assume that the minimum-phase invariant zeros are all distinct (i.e., their algebraic multiplicity is one):

\begin{assumption}
 \label{Ass2}
System $\Sigma$ has no coincident minimum-phase invariant zeros.
\end{assumption}

This assumption does not lead to a significant loss of generality. In fact, the case of coincident zeros can be dealt with by using a procedure that is similar in spirit to that outlined in \cite{Schmid-NNP-13}. 

\begin{lemma}{\sc (\cite{SICON,CDC})}.
\label{cor2}
Let $r = \dim \gR^\star$ and assume that Assumption~\ref{Ass2} holds. Let $\gZ_g=\{z_{r+1},z_{r+2},\ldots,z_{r+t}\}$ be the $s_z$-conformably ordered set of minimum-phase invariant zeros of $\Sigma$. Let $\gL = \{\mu_1,\ldots,\mu_r\}$ be $s$-conformably ordered such that $\gL \cap \gZ_g=\emptyset$. Let $K = \diag \{k_1,\ldots,k_{r}\}$ be defined as in Lemma~\ref{Rs}.
Let $H \defi \diag \{h_{r+1},\ldots,h_{r+t}\}$, where $h_i \in \complex^{\dim \left(\ker P_{\scriptscriptstyle \Sigma}(z_i)\right)}$ for all {$i \in \{r+1,\ldots,r+t\}$ and $\overline{h}_i=h_{i+1}$ for all odd $i-r \in \{1,\ldots,2\,s_z-1\}$.}
Let 
\beann
{M}_{\scriptscriptstyle K,H} = \bmat{c|c|c|c|c|c|c} \!\! N_{\scriptscriptstyle \Sigma}(\mu_1) & \ldots & N_{\scriptscriptstyle \Sigma}(\mu_r) & N_{\scriptscriptstyle \Sigma}(z_{r+1}) & N_{\scriptscriptstyle \Sigma}(z_{r+2}) & \ldots & N_{\scriptscriptstyle \Sigma}(z_{r+t}) \!\! \emat \diag\{K,H\}
\eeann
and let for all $i \in \{1,\ldots,r+t\}$
\beann
\bmat{c} {v}_{\scriptscriptstyle K,H,i} \\ {w}_{\scriptscriptstyle K,H,i} \emat=\left\{ \begin{array}{ll}
\!\! \mathfrak{Re} \{M_{\scriptscriptstyle K,H}^i\} & \textrm{if $i-r \in \{1,\ldots,2\,s_z\}$ is odd} \\[1mm]
\!\! \mathfrak{Im} \{M_{\scriptscriptstyle K,H}^i\} & \textrm{if
 $i-r \in \{1,\ldots,2\,s_z\}$ is even} \\[1mm]
M_{\scriptscriptstyle K,H}^i & \textrm{if $i \in \{1,\ldots,r\}\cup \{r+2\,s_z+1,\ldots,r+t\}$} \end{array} \right.
\eeann
where ${v}_{\scriptscriptstyle K,H,i}$ are $n$-dimensional and ${w}_{\scriptscriptstyle K,H,i}$ are $m$-dimensional for all $i$.
Finally, let
\bea
\label{tXKYK}
V_{\scriptscriptstyle K,H} \ns&\ns = \ns&\ns [\, v_{\scriptscriptstyle K,H,1} \;\; \ldots \;\; v_{\scriptscriptstyle K,H,r} \;\;{v}_{\scriptscriptstyle K,H,r+1} \;\; \ldots \;\;{v}_{\scriptscriptstyle K,H,r+t} \,], \label{XK1} \\
W_{\scriptscriptstyle K,H} \ns&\ns = \ns&\ns [\, w_{\scriptscriptstyle K,H,1} \;\; \ldots \;\; w_{\scriptscriptstyle K,H,r} \;\; {w}_{\scriptscriptstyle K,H,r+1} \;\; \ldots \;\; {w}_{\scriptscriptstyle K,H,r+t} \,]. \label{YK1}
\eea
Thus, {\bf (i)} For almost every choice of  $K = \diag \{k_1,\ldots,k_r\}$ and $H = \diag \{h_{r+1},\ldots,h_{r+t}\}$ we have $\rank V_{\scriptscriptstyle K,H}=r+t$; {\bf (ii)} If $K$ and $H$ are such that $\rank V_{\scriptscriptstyle K,H}=r+t$, the matrix $V_{\scriptscriptstyle K,H}$ is a basis matrix for $\gV^\star_g$ adapted to $\gR^\star$; {\bf (iii)} 
The set of all friends of $\gV^\star_g$ such that $\sigma(A+B\,F\,|\,\gV^\star_g)=\gL \cup \gZ_g$ is parameterised in $K$ and $H$ as $F_{\scriptscriptstyle K,H}=W_{\scriptscriptstyle K,H}\,V_{\scriptscriptstyle K,H}^\dagger$, 
where ${K,H}$ are such that $\rank V_{\scriptscriptstyle K,H}=r+t$.\\[-6mm]
\end{lemma}

\begin{remark}
\label{agg1}
We now show how we can build a basis for $\gV^\star_g$ using the result of Lemma \ref{cor2}. Let $\{\mu_1,\mu_2,\ldots,\mu_{g}\}\subset \complex_g$ include the minimum-phase invariant zeros. Let this set be $s$-conformably ordered. 
In view of Lemma \ref{cor2}, for all $i\in \{1,\ldots,g\}$ we can find vectors $\bsmat \phi_i \\[1mm] \psi_i \esmat$ such that $\bsmat A-\mu_i \,I_n & B \\[1mm] C & D \esmat\bsmat \phi_i \\[1mm] \psi_i \esmat=0$ for all $i \in \{1,\ldots,g\}$ 
such that, by defining $v_{i}=\mathfrak{Re}\{ \phi_i\}$ when $i \le 2\,s$ is odd, $v_{i}=\mathfrak{Im}\{ \phi_i\}$ when $i \le 2\,s$ is even, and $v_{i}= \phi_i$ when $i \in \{2\,s+1,\ldots,g\}$, the matrix $[\,v_1 \;\ldots\; v_{g}\,]$ is of full column-rank.
We can define $w_{1},\ldots,w_g$ similarly using $\psi_i$ instead of $\phi_i$. 
Then, by virtue of Lemma \ref{cor2} we have $\spanR \{v_{1},\ldots,v_g\}=\gV^\star_g$.
\end{remark}
\begin{corollary}
\label{agg3}
Let $r = \dim \gR^\star$ and assume that Assumption~\ref{Ass2} holds. Let $\gZ_g=\{z_{1},\ldots,z_{t}\}$ be the $s_z$-conformably ordered set of minimum-phase invariant zeros of $\Sigma$. When $\mu\in \complex$, let us define 
$\gR^\star (\mu) \stackrel{\text{\tiny def}}{=} \left\{ \spanR \{\mathfrak{Re}\{v\},\mathfrak{Im}\{v\}\} \subseteq \gX \,\Big|\,\, \exists w \in \complex^m \;:\;\;
\left[ \begin{smallmatrix} A-\mu\,I_n & B \\[1mm] C & D \end{smallmatrix} \right]\left[ \begin{smallmatrix} v \\[1mm] w \end{smallmatrix} \right]=0 \right\}$.
There holds
\[
\gV^\star_g=\gR^\star+\gR^\star(z_{1})+\ldots+\gR^\star(z_{t}).
\]
\end{corollary}

\section{Solution to {Problem~\ref{prob:mono2}}}

Our aim in this section is to provide tractable and constructive necessary and sufficient conditions for the existence of a solution to the problem of global monotonicity. 

\subsection{A first necessary and sufficient condition for {Problem~\ref{prob:mono2}}}
As explained above, in order to achieve a globally monotonic step response we need to find the feedback matrix $F$ that evenly distributes $p$ of the $n$ closed-loop modes into the $p$ components of the tracking error, and renders the remaining $n-p$ modes invisible at the tracking error. The number of closed-loop modes that can be made invisible by state feedback equals the dimension of the subspace $\gV^\star_g$. 
Thus, for the tracking control problem with global monotonicity to be solvable we need the condition $\dim \gV^\star_g \ge n-p$ to be satisfied. This condition is only necessary, because we need also the {linearly independent} vectors $v_1,\ldots,v_p$ obtained with the procedure indicated above to be linearly independent of $\gV^\star_g$. In the case in which $\dim \gV ^\star_g > n-p$ holds, if it is possible to find {linearly independent} vectors $v_1,\ldots,v_p$ that are independent of $\gV^\star_g$, then not only is the monotonic tracking control problem solvable, but we are potentially able to also obtain a response that achieves instantaneous tracking in some outputs.


We now make a simplifying technical assumption, which in view of the discussion above amounts to putting ourselves in a ``worse-case scenario'' of all the possible situations in which the tracking problem is solvable. This assumption is made for the sake of simplicity:
\begin{assumption}
\label{assVg}
$\dim \gV^\star_g = n-p$.
\end{assumption}


Let $\l \in \real$. For all $j \in \{1,\ldots,p\}$ we define 
\begin{equation}
\label{Rhat}
	\hat{\gR}_j (\lambda) \stackrel{\text{\tiny def}}{=} \left\{ v \in \gX \,\Big|\,\, \exists \beta \in \mathbb{R}\setminus \{0\}, \, \exists w \in \gU \;:\;\;
\left[ \begin{smallmatrix} A-\lambda\,I_n & B \\[1mm] C & D \end{smallmatrix} \right]\left[ \begin{smallmatrix} v \\[1mm] w \end{smallmatrix} \right]=\left[ \begin{smallmatrix} 0 \\[1mm] \beta e_j \end{smallmatrix} \right] \right\}.
\end{equation}

It is easy to see that, given $\lambda \in \real$, the set $\hat{\gR}_j (\lambda)$ is not a subspace of $\gX$.


The following lemma provides a necessary and sufficient condition for Problem~\ref{prob:mono2} to admit solutions in terms of the sets $\hat{\gR}_j (\lambda)$ defined above.

\begin{lemma}
\label{lem:Rosej}
{Let $\mathfrak{L}_p =(\lambda_1, \lambda_2, \cdots, \lambda_p) \in \real^p$ be such that in the continuous time $\lambda_j \in \real^-$, and in the discrete time $\lambda_j \in (0,1)$.} 
Under Assumptions~\ref{Ass1} and \ref{assVg}, Problem~\ref{prob:mono2} admits solution if and only if there exist $(v_{1},\cdots,v_{p})\in\hat{\gR}_{1} (\lambda_{1})\times\cdots\times \hat{\gR}_{p} (\lambda_{p})$ satisfying
\begin{equation} 
\label{eq:Rosej}
 \gV^\star_g + \spanR \{ v_{1},\cdots,v_{p}\} = \gX.
\end{equation}
\end{lemma}
\proof
Let us consider for the sake of argument the continuous time. The discrete case follows with the obvious substitutions. First, we show sufficiency. Since we are assuming $\dim \gV^\star_g=n-p$ and that (\ref{eq:Rosej}) holds, then $\dim (\spanR \{v_{1},\cdots,v_{p}\})=p$, which means that $\{v_{1},\cdots,v_{p}\}$ are linearly independent. From (\ref{Rhat}), {there exists $w_i\in \gU$ such that} 
$\bsmat A-\lambda_i \,I_n & B \\[1mm] C & D \esmat \bsmat v_i \\[1mm] w_i \esmat= \bsmat 0 \\[1mm] \beta_i\,e_i \esmat$ for $i \in \{1,\ldots,p\}$ where $\beta_i \neq 0$. {{We now build a basis for $\gV^\star_g$ as shown in Corollaries \ref{agg} and \ref{agg3}. Since $\dim \gV^\star_g=n-p$,} let $\{\mu_1,\mu_2,\ldots,\mu_{n-p}\}\subset \complex_g$ include the minimum-phase invariant zeros. Let this set be $s$-conformably ordered. 
Using the consideration in Remark \ref{agg1}, we find $\spanR \{v_{p+1},\ldots,v_n\}=\gV^\star_g$, so that from (\ref{eq:Rosej}) the set $\{v_1,\ldots,v_p,v_{p+1},\ldots,v_n\}$ is linearly independent.} Thus, constructing $\{w_{p+1},\ldots,w_n\}$ also as in Remark \ref{agg1}, the feedback matrix $F \defi [\,w_1\,\ldots\,w_n\,]\,[\,v_1\,\ldots\,v_n\,]^{-1}$ satisfies
\beann
(A+B\,F) \!\!\! \!\!\! \!\!\! \!\!\! \ns&\ns \!\!\! \!\!\! \ns&\ns[\,v_1\;\; \ldots\;\; v_p\;|\;v_{p+1}\;\;v_{p+2}\;|\;\ldots\;|\;v_{p+2\,s-1}\;\;v_{p+2\,s}\;|\;v_{p+2\,s+1} \;\; \ldots\;\; v_n\,] \\
 \ns&\ns = \ns&\ns
 \diag \left\{ \lambda_1,\ldots, \lambda_p,
 \bsmat \mathfrak{Re}\{\mu_1\} & -\mathfrak{Im}\{\mu_1\} \\[1mm]
 \mathfrak{Im}\{\mu_1\} & \mathfrak{Re}\{\mu_1\} \esmat, \ldots, 
 \bsmat \mathfrak{Re}\{\mu_{2\,s-1}\} & -\mathfrak{Im}\{\mu_{2\,s-1}\} \\[1mm]
 \mathfrak{Im}\{\mu_{2\,s-1}\} & \mathfrak{Re}\{\mu_{2\,s-1}\} \esmat, \mu_{2\,s+1},\ldots,\mu_{n-p}\right\}\\ 
(C+D\,F)\,v_i \ns&\ns = \ns&\ns \left\{ \begin{array}{ll} \beta_i\,e_i & i \in \{1,\ldots,p\} \\
0 & i \in \{p+1,\ldots,n\}
\end{array} \right.
\eeann
 Let $\xi_0=\xi(0)$ be the initial error state, and define $\alpha \defi [\,v_1\,\ldots\,v_n\,]^{-1}\,\xi_0$. We find
\beann
\epsilon(t) \ns&\ns = \ns&\ns (C+D\,F)\,\exp\left[(A+B\,F)\,t \right]\,\xi_0 \\
\ns&\ns = \ns&\ns \sum_{i=1}^p \exp(\lambda_i\,t)\,(C+D\,F)\,v_i\,\alpha_i = \sum_{i=1}^p \beta_i\,e_i\,\exp(\lambda_i\,t)\,\alpha_i=\bsmat \beta_1\,\alpha_1\,\exp(\lambda_1\,t) \\[-1mm] \vdots \\[1mm] \beta_p\,\alpha_p\,\exp(\lambda_p\,t)\esmat,
 \eeann
 so that each component of $\epsilon(t)$ is given by a single exponential, and is therefore monotonic.\\
 {
Let us now consider necessity. If Problem~\ref{prob:mono2} admits the solution $F$, by Lemma \ref{lem:mono}, the tracking error has a single closed-loop mode per component, i.e., it is in the form given by (\ref{risp}). This implies that the remaining $n-p$ closed-loop modes (which are asymptotically stable because $F$ is stabilising) must disappear from the tracking error. Hence, $\dim \gV^\star_g \ge n-p$ (and in particular $\dim \gV^\star_g = n-p$ in view of Assumption \ref{assVg}). Let us define, as in Lemma \ref{lem:mono}, by $V \defi [\,v_1\;\;v_2\;\;\ldots\;\;v_n\,]$ the eigenvector matrix of $A+B\,F$, so that $V$ is invertible (recall that the closed-loop eigenvalues are distinct as established in the proof of Lemma \ref{lem:mono}). No generality is lost by assuming that $\{v_{p+1},\ldots,v_n\}$ is a basis for $\gV^\star_g$. Let us define $\xi_0 \in \gX$ and $\alpha =V^{-1} \xi_0$ as in Lemma \ref{lem:mono}. Then
\beann
\epsilon(t)=(C+D\,F)\,[\,v_1\,e^{\l_1\,t}\,\alpha_1+\ldots+v_p\,e^{\l_p\,t}\,\alpha_p\,]=\bsmat \gamma_1\,\exp(\lambda_1\,t) \\[-1mm] \vdots \\[1mm] \gamma_p\,\exp(\lambda_p\,t)\esmat
\eeann
for some $\gamma_1,\ldots,\gamma_p\in \real$, since $(C+D\,F)\,v_i=\{0\}$ for all $i \in \{p+1,\ldots,n\}$. 
Consider $\epsilon_j(t)=[\,0\;\;\ldots\;\; \gamma_j\,\exp(\lambda_j\,t)\;\; \ldots\;\; 0\,]^\top$, which corresponds to the (non-unique) initial error state $\xi_{0,j}=V\,e_j$, where $e_j$ is the $j$-th canonical basis vector of $\gX$. It is easy to see that $\xi_{0,j}$ spans an output-nulling subspace of the system {$(A,B,C_{(j)},D_{(j)})$} obtained by removing the $j$-th output because all components of $\epsilon_j(t)$ except the $j$-th are zero. Since in $\epsilon_j(t)$ there is only one mode, we have $(A+B\,F)\,\xi_{0,j}=\lambda_j\,\xi_{0,j}$ and $(C_{(j)}+D_{(j)}\,F)\,\xi_{0,j}=0$. The latter implies $(C+D\,F)\,\xi_{0,j}=\gamma_j\,e_j$ for a certain $\gamma_j$, which cannot be zero, because this would imply $\dim \gV^\star_g > n-p$ against Assumption \ref{assVg}. Thus, we may define $v_j=\frac{1}{\gamma_j}\,\xi_{0,j}$ and $w_j=\frac{1}{\gamma_j}\,F\,\xi_{0,j}$, which satisfy (\ref{ML0}). By superposition, we need  $\{v_1,\ldots,v_p\}$ to span a $p$-dimensional subspace of $\gX$ independent of $\gV_g^\star$.
}
\endproof


\begin{remark}
Whenever (\ref{eq:Rosej}) is satisfied, Problem~\ref{prob:mono2} can be solved with an arbitrary convergence rate. At first glance, this property seems to be in contrast with the fact that the pair $(A,B)$ has not been assumed to be completely reachable, but only stabilisable. In other words, one may argue that the uncontrollable modes (which are asymptotically stable), may limit the convergence rate. However, it is easy to see that this is not the case. Indeed, from the right invertibility of the quadruple $(A,B,C,D)$, one can conclude that every uncontrollable
eigenvalue of the pair $(A,B)$ is also an invariant zero of $\Sigma$. \footnote{This can be seen by observing that an uncontrollable eigenvalue $\lambda$ of $(A,B)$ either belongs to
 $\sigma (A+B\,\Phi \,|\,\gX/\gV^\star\!+\!\gR_0)$ or to $\sigma (A+B\,\Phi \,|\,\gV^\star / \gR^\star)$, where $\gR_0=\langle A, \ima B \rangle$ is the reachability subspace of the pair $(A,B)$, i.e., the smallest $A$-invariant subspace containing the range of $B$, and $\Phi$ is any
 friend of $\gV^\star$. Since $\gR_0$ is contained in the smallest input-containing subspace $\gS^\star$ of $\Sigma$ \cite[Chapter 8]{Trentelman-SH-01}, and the right-invertibility is equivalent to the condition $\gV^\star + \gS^\star=\gX$ since the matrix $[\,C\;\;\;D\,]$ has been assumed to be of full row-rank \cite[Theorem 8.27]{Trentelman-SH-01}, we also have $\gV^\star + \gR_0=\gX$. Hence, $\lambda \in \sigma (A+B\,\Phi \,|\,\gV^\star / \gR^\star)$, i.e., $\lambda\in \gZ$. }
Hence, every uncontrollable eigenvalue of the pair $(A,B)$ is rendered invisible at the tracking error, and therefore it does not limit the rate of convergence. 
It is also worth observing that  there is freedom in the choice of the closed-loop eigenvalues associated with $\gR^\star$, when computing a basis matrix for $\gV^\star_g$. Even though these eigenvalues are invisible at the tracking error (and hence any choice will be correct as long as they are asymptotically stable and distinct from the minimum-phase invariant zeros) this freedom may be important for the designer, since the selection of closed-loop eigenvalues affects other considerations like control amplitude/energy. Thus, it is worth emphasising that the designer has complete freedom to chose any set of stable eigenvalues provided the minimum-phase invariant zeros are included, and provided at least $p$ of these meet the desired convergence rate.
\\[-5mm]
\end{remark}

Lemma \ref{lem:Rosej} already provides a set of necessary and sufficient conditions for the solvability of the globally monotonic tracking control problem. However, such conditions are not easy to test, because they are given in terms of the sets $\hat{\gR}_j(\lambda_j)$ which are not, in general, subspaces of $\gX$.
The tools that we now present are aimed at replacing $\hat{\gR}_j(\lambda_j)$ in condition (\ref{eq:Rosej}) with particular reachability subspaces of the state-space, which we now define.
{As in the proof of Lemma \ref{lem:Rosej}, 
for each output} $j\in \{1,\ldots,p\}$ we introduce $\Sigma_j=(A,B,C_{(j)},D_{(j)})$ as the quadruple in which $C_{(j)}\in \real^{(p-1) \times n}$ and $D_{(j)}\in \real^{(p-1) \times m}$ are obtained by eliminating the $j$-th row from $C$ and $D$, respectively. 
We observe that the right invertibility of the quadruple $(A,B,C,D)$ guarantees that the set of invariant zeros of $\Sigma$ contains the set of invariant zeros of $\Sigma_j$ for any $j\in \{1,\ldots,p\}$.
The largest output nulling reachability subspace of $\Sigma_j$ is denoted by $\gR^\star_j$. {Similarly to what was done for $\gR^\star$ in Corollary \ref{agg}, {for any distinct set $\{\mu_{1},\cdots,\mu_{r_j}\} \subset \real \setminus \gZ$,} we decompose $\gR^\star_j$ as}
\begin{equation} \label{eq:gRdec}
\gR^{\star}_j=\gR^{\star}_j(\mu_{1})+\cdots+\gR^{\star}_j(\mu_{r_j}),
\end{equation}
where $r_j = \dim \gR^\star_j$ {and}
\begin{equation} \label{eq:Rsj}
	\gR^\star_j (\mu_i) \stackrel{\text{\tiny def}}{=} \left\{ v \in \gX \,\Big|\,\, \exists w \in \gU \;:\;\;
\left[ \begin{smallmatrix} A-\mu_i\,I_n & B \\[1mm] C_{(j)} & D_{(j)} \end{smallmatrix} \right]\left[ \begin{smallmatrix} v \\[1mm] w \end{smallmatrix} \right]=0 \right\}.
\end{equation}

\begin{remark}
{As established in Corollary \ref{agg}, a spanning set} for $\gR^\star_j (\mu_i)$ is given by the columns of $V_i$, where $V_i$ is the upper part of a basis matrix $\bsmat V_i \\[1mm] W_i \esmat$ of $\left[ \begin{smallmatrix} A-\mu_i\,I_n & B \\[1mm] C_{(j)} & D_{(j)} \end{smallmatrix} \right]$. However, differently from $\gR^\star(\mu)$, this time it is not guaranteed that $V_i$ obtained in this way is of full column-rank, because the matrix $\left[ \begin{smallmatrix} B \\[1mm] D_{(j)} \end{smallmatrix} \right]$ may very well have a non-trivial kernel. 
\end{remark}

The relationship between $\hat{\gR}_j (\mu)$ and $\gR^\star_j(\mu)$ is stated through the two following results.

\begin{proposition} 
\label{prop:RhRs}
Let $\mu \in \real \setminus \gZ$. 
For all $j \in \{ 1,\cdots , p\}$, there holds
\begin{equation}
	\gR^\star_j (\mu) = \hat{\gR}_j (\mu) \cup \gR^\star (\mu)
\end{equation}
\end{proposition}

\proof
First, we prove that $\gR^\star_j (\mu) \supseteq \hat{\gR}_j (\mu) \cup \gR^\star (\mu)$. To this end, we first show that $\hat{\gR}_j (\mu) \subseteq \gR^\star_j (\mu)$. Let $v \in \hat{\gR}_j (\mu)$. There exist $w \in \gU$ and $\beta \in \mathbb{R}\setminus \{0\}$ such that
$\left[ \begin{smallmatrix} A-\mu\,I_n & B \\[1mm] C & D \end{smallmatrix} \right]\left[ \begin{smallmatrix} v \\[1mm] w \end{smallmatrix} \right]=\left[ \begin{smallmatrix} 0 \\[1mm] \beta e_j \end{smallmatrix} \right]
$,
which implies in particular that $C_{(j)}\,v+D_{(j)}\,w=0$. Hence, $v \in \gR^\star_j (\mu)$.
{We now show that $\gR^\star (\mu) \subseteq \gR^\star_j (\mu)$. Let $v \in \gR^\star (\mu)$. Then, there exist $w \in \gU$ such that
$\left[ \begin{smallmatrix} A-\mu\,I_n & B \\[1mm] C & D \end{smallmatrix} \right]\left[ \begin{smallmatrix} v \\[1mm] w \end{smallmatrix} \right]=\left[ \begin{smallmatrix} 0 \\[1mm] 0 \end{smallmatrix} \right]$,
which again implies that $C_{(j)}\,v+D_{(j)}\,w=0$, so that $v \in \gR^\star_j (\mu)$. Hence, $\hat{\gR}_j (\mu) \cup \gR^\star (\mu) \subseteq \gR^\star_j (\mu)$ holds.}
We now show that $\gR^\star_j (\mu) \subseteq \hat{\gR}_j (\mu) \cup \gR^\star (\mu)$. Let
 $v$ be an element of $\gR_j^\star (\mu)$. Then, there exists a $w \in \gU$ such that
$\left[ \begin{smallmatrix} A-\mu\,I_n & B \\[1mm] C_{(j)} & D_{(j)} \end{smallmatrix} \right]\left[ \begin{smallmatrix} v \\[1mm] w \end{smallmatrix} \right]=0$. 
Let $\beta = C_{j}\,v+D_{j}\,w$. Then, $\left[ \begin{smallmatrix} A-\mu\,I_n & B \\[1mm] C & D \end{smallmatrix} \right]\left[ \begin{smallmatrix} v \\[1mm] w \end{smallmatrix} \right]=\left[ \begin{smallmatrix} 0 \\[1mm] \beta e_j \end{smallmatrix} \right]$. If $\beta \neq 0$, we have $v \in \hat{\gR}_j (\mu)$, whereas if $\beta=0$, we find $v \in \gR^\star (\mu)$. Thus, $v \in \hat{\gR}_j (\mu) \cup \gR^\star (\mu)$.
\endproof

{
\begin{proposition} \label{prop:meashat}
Let $\mu \in \real \setminus \gZ$. For all $j \in \{ 1,\cdots , p\}$, there holds
\begin{equation} \label{eq:meashat}
{m}_{\scriptscriptstyle \gR^\star_j (\mu)} (\gR^\star_j (\mu) \setminus \hat{\gR}_j(\mu)) =0. 
\end{equation}
\end{proposition}
}

\proof
Since $\Sigma$ is right invertible and $\mu$ is not an invariant zero, the inclusion $\gR^\star (\mu) \subseteq \gR_j^\star (\mu)$ deriving from Proposition \ref{prop:RhRs} becomes $\gR^\star (\mu) \subset \gR_j^\star (\mu)$. Indeed, in such a case, $[\, C_{j}\;\;\; D_{j}\,]$ is linearly independent from every row of $\left[ \begin{smallmatrix} A-\mu\,I_n & B \\[1mm] C_{(j)} & D_{(j)} \end{smallmatrix} \right]$. This implies that $\dim \gR^\star (\mu) < \dim \gR_j^\star (\mu)$, so that ${m}_{\scriptscriptstyle \gR^\star_j (\mu)} (\gR^\star(\mu)) =0$. 
Moreover, 
Proposition \ref{prop:RhRs} ensures 
 that $\gR^\star_j (\mu) \setminus \hat{\gR}_j (\mu) \subseteq \gR^\star (\mu)$, which in general does not hold as an equality since $\gR^\star (\mu)$ and $\hat{\gR}_j (\mu)$ may very well have non-zero intersection. Thus, \eqref{eq:meashat} follows readily.
\endproof





 Roughly speaking, this result, {together with Proposition~\ref{prop:RhRs}}, implies that $\gR^\star_j (\mu)$ is coincident with $\hat{\gR}_j (\mu)$ modulo a set of points belonging to a proper algebraic variety within $\gR^\star_j (\mu)$. 
 This essential step justifies the fact that from now on we will use 
 $\gR^\star_j (\mu)$, instead of $\hat{\gR}_j (\mu)$, to establish constructive necessary and sufficient condition for our tracking problem.

\subsection{A tractable condition for the solution of Problem~\ref{prob:mono2}}

{
Let $\hat{\mathfrak{L}}_p$ be the set of all $p$-tuples $(\lambda_{1},\cdots,\lambda_{p})\in\mathbb{R}^{p}$ such that for all $i \in \{1,\cdots,p\}$ we have $\lambda_i\notin {\cal Z}$, and $\lambda_i \in \real^-$ or $\lambda_i \in (0,1)$ in the continuous or in the discrete time, respectively.
}

\begin{theorem} \label{th:kchoice}
Let ${\frak L}_p = (\lambda_1,\ldots,\lambda_p) \in \hat{\mathfrak{L}}_p$. 
Problem~\ref{prob:mono2} admits solution if and only if
\begin{equation}
\forall\,{S}\in2^{\{1,\ldots,p\}},\qquad\dim\left(\gV_{g}^{\star}+\sum_{j\in{S}}\gR_{j}^{\star}(\lambda_{j})\right)\ge n-p+{\rm {card}({S})}. \label{eq:condL}
\end{equation}
\end{theorem}

\proof
We begin by defining the propositions
\begin{eqnarray}
&&\mathfrak{q}_1: \qquad
 \gV^\star_g + \spanR \{ v_{1},\cdots,v_{p}\} = \gX, \\
&& \mathfrak{q}_2: \qquad ( v_{1},\cdots,v_{p} )\in \hat{\gR}_{1} (\lambda_{1})\times\cdots\times \hat{\gR}_{p} (\lambda_{p}), \label{eq:condLoc} 
\end{eqnarray}
and the sets
\begin{eqnarray}
\gT_{1} & \defi & \left\{ ( v_{1},\cdots,v_{p} )\in \gR_{1}^{\star}(\lambda_{1}) \times \cdots \times \gR_{p}^{\star}(\lambda_{p}) \big| \neg \mathfrak{q}_1 \right\}, \nn \\
\gT_{2} & \defi & \left\{ ( v_{1},\cdots,v_{p} )\in \gR_{1}^{\star}(\lambda_{1}) \times \cdots \times \gR_{p}^{\star}(\lambda_{p}) \big| \neg \mathfrak{q}_2 \right\}, \nn \\
\gT_{1,2} & \defi & \left\{ ( v_{1},\cdots,v_{p} )\in \gR_{1}^{\star}(\lambda_{1}) \times \cdots \times \gR_{p}^{\star}(\lambda_{p}) \big| \neg (\mathfrak{q}_1 \land \mathfrak{q}_2) \right\}. \nn 
\end{eqnarray}

Suppose that (\ref{eq:condL}) is satisfied. We define for the sake of conciseness $\tilde{\gR} \defi \gR_{1}^{\star}(\lambda_{1}) \times \cdots \times \gR_{p}^{\star}(\lambda_{p})$. Then, 
${m}_{\scriptscriptstyle \tilde{\gR}}(\gT_{1})=0$ and ${m}_{\scriptscriptstyle \tilde{\gR}}(\gT_{2})=0$ are ensured by Lemma~\ref{lem:kchoice} in Appendix~A, and Proposition~\ref{prop:meashat}. Since $\gT_{1,2} = \gT_{1} \cup \gT_{2}$, it follows that ${m}_{\scriptscriptstyle \tilde{\gR}} (\gT_{1,2}) \leq {m}_{\scriptscriptstyle \tilde{\gR}} (\gT_{1}) + {m}_{\scriptscriptstyle \tilde{\gR}} (\gT_{2})=0$. This is equivalent to saying that for almost all $( v_{1},\cdots,v_{p} )\in \gR_{1}^{\star}(\lambda_{1}) \times \cdots \times \gR_{p}^{\star}(\lambda_{p})$ both $\mathfrak{q}_1$ and $\mathfrak{q}_2$ hold. According to Lemma~\ref{lem:Rosej}, this proves that Problem~\ref{prob:mono2} admits solution.
Suppose now that (\ref{eq:condL}) is not satisfied and note that Proposition~\ref{prop:RhRs} ensures that $\hat{\gR}_{1} (\lambda_{1})\times\cdots\times \hat{\gR}_{p} (\lambda_{p})\subseteq \gR_{1}^{\star}(\lambda_{1}) \times \cdots \times \gR_{p}^{\star}(\lambda_{p})$ since every $\lambda_j$ belongs to $\real \setminus \gZ$.
In such a case, the second statement of Lemma~\ref{lem:kchoice} guarantees that there is no $( v_{1},\cdots,v_{p} )$ verifying $\mathfrak{q}_1$ which belongs to $\gR_{1}^{\star}(\lambda_{1}) \times \cdots \times \gR_{p}^{\star}(\lambda_{p})\supseteq\hat{\gR}_{1} (\lambda_{1})\times\cdots\times \hat{\gR}_{p} (\lambda_{p})$. Thus, Problem~\ref{prob:mono2} does not admit solution in view of Lemma~\ref{lem:Rosej}.
\endproof

\begin{example}
Consider the system in Example \ref{exe0}. If we denote by $\{e_1,\ldots,e_5\}$ the canonical basis in $\gX=\real^5$, it is easily verified that $\gR^\star_1=\gR^\star_3= \spanR\{e_2,e_3,e_4,e_5\}$, $\gR^\star_2= \spanR\{e_3,e_4,e_5\}$. We recall that $\gV^\star_g=\ima \bsmat -2 && 2/3 && -41/22 && 0 && -1/11 \\[1mm] 0 && 0 && 0 && 1 && 0 \esmat^\top$. In this case, (\ref{eq:condL}) can be written as:
\beann
\dim (\gV^\star_g+\gR^\star_j) \ns&\ns \ge \ns&\ns n-p+1 \quad \forall \,j \in \{1,2,3\} \\
\dim (\gV^\star_g+\gR^\star_i+\gR^\star_j) \ns&\ns \ge \ns&\ns n-p+2 \quad \forall \,i,j \in \{1,2,3\}\, \text{such that $i \neq j$} \\
\dim (\gV^\star_g+\gR^\star_1+\gR^\star_2+\gR^\star_3) \ns&\ns \ge \ns&\ns n-p+3.
\eeann
In the present case, these conditions are verified. Indeed, we find $\gV^\star_g+\gR^\star_1=\gV^\star_g+\gR^\star_3=\gX$, the dimension of $\gV^\star_g+\gR^\star_2$ is $4$, and $\gV^\star_g+\gR^\star_2+\gR^\star_3=\gX$.
\end{example}

\subsection{Computation of the gain feedback} \label{sec:algo}
Let $V_g$ be a basis matrix for $\gV^\star_g$.
We first consider the case in which $V_g$ has $h=n-p$ columns. Let
$V_g=\left[ v_{g,1}\;\;\,v_{g,2}\;\;\ldots \;\;v_{g,h}\right]$ and 
$W_g=\left[ w_{g,1}\;\;\,w_{g,2}\;\;\ldots \;\;w_{g,h}\right]$, which satisfy
$\bsmat A-\mu_i\,I_n & B \\[1mm] C & D \esmat \bsmat v_{g,i} \\[1mm] w_{g,i} \esmat=0$ (here we assume for the sake of simplicity that all the $\mu_i$ are real, but in the case of complex conjugate minimum-phase invariant zeros, one can apply the construction of Corollaries \ref{agg} and \ref{agg3} with the obvious modifications).
The necessary and sufficient condition $\gV^\star_g + \spanR \{ v_{1},\cdots,v_{p}\} = \gX$ is satisfied with {$(v_{1},\ldots,v_p)\in\gR_{1}^{\star}(\lambda_{1}) \times \cdots \times \gR_{p}^{\star}(\lambda_{p})$ and, hence, there exists} $\{w_1,\ldots,w_p\}$ such that $\bsmat A-\lambda_i\,I_n & B \\[1mm] C_{(i)} & D_{(i)} \esmat \bsmat v_{i} \\[1mm] w_{i} \esmat=0$.
Since this condition is equivalent to the condition $\rank [\,v_1\;\;\ldots\;\;v_p\;\;V_g\,]=n$, we can compute $F=W\,V^{-1}$. 
{From Proposition~\ref{prop:meashat}, $(v_{1},\ldots,v_p)\in\hat{\gR}_{1} (\lambda_{1})\times\cdots\times \hat{\gR}_{p} (\lambda_{p})$ generically holds.
} Then, $F$ ensures that $(A+B\,F)\,v_{g,i}=\mu_i\,v_{g,i}$ and $(C+D\,F)\,v_{g,i}=0$ for all $i \in \{1,\ldots,n-p\}$, 
and that there exist $\beta_i \neq 0$ such that 
$\bsmat A-\lambda_i\,I_n & B \\[1mm] C & D \esmat \bsmat v_{i} \\[1mm] w_{i} \esmat=\bsmat 0 \\ \beta_{i} \,e_i\esmat$, which in turn gives $(A+B\,F)\,v_{i}=\lambda_i\,v_{i}$ and $(C+D\,F)\,v_{i}=\beta_i\,e_i$ for all $i \in \{1,\ldots,p\}$. Therefore, $\sigma(A+B\,F)=\{\lambda_1,\ldots,\lambda_p,\mu_1,\ldots,\mu_{n-p}\}$ and 
\[
\epsilon(t)=(C+D\,F)\,\exp({\lambda_1\,t})\,v_1\,\gamma_1+\ldots+(C+D\,F)\,\exp({\lambda_p\,t})\,v_p\,\gamma_p=\bmat{c} \beta_1\,\gamma_1\,\exp({\lambda_1\,t}) \\ \vdots \\ 
\beta_p\,\gamma_p\,\exp({\lambda_p\,t}) \emat
\]
{for some $\gamma_1,\ldots,\gamma_p \in \real$ as required.} We now consider the case where Assumption \ref{assVg} does not hold. In other words, since we know that $\dim \gV^\star_g\ge n-p$ is a necessary solvability condition, we now assume $\dim \gV^\star_g=h>n-p$. 

\begin{proposition}
 \label{th:kchoiceGen} 
Let ${\frak L}_p = (\lambda_1,\ldots,\lambda_p) \in \hat{\mathfrak{L}}_p$. Problem~\ref{prob:mono2} admits solution if and only if there exists {a set $\delta \subseteq\{1,\cdots,p\}$ satisfying $\mathrm{card}\, \delta = n-h$ and}
\begin{equation} \label{eq:condLh}
\forall\,{S}\in2^{\delta} , \quad \dim\left(\gV_{g}^{\star}+\sum_{j\in S}\gR_{j}^{\star}(\lambda_{j})\right)\ge h+{\rm {card}({S})}.
\end{equation}
\end{proposition}

\proof
Using the same argument of the proof of Lemma~\ref{lem:Rosej},  Problem~\ref{prob:mono2} is seen to admit solutions if and only if there exists {$\delta \subseteq\{1,\cdots,p\}$ with $\mathrm{card}\, \delta = n-h$ and a bijective map $\beta: \{1,\cdots,n-h\}\longrightarrow \delta$ } such that $(v_{\beta(1)},\cdots,v_{\beta(n-h)})\in\hat{\gR}_{\beta(1)} (\lambda_{\beta(1)})\times\cdots\times \hat{\gR}_{\beta(n-h)} (\lambda_{\beta(n-h)})$ satisfies $\gV^\star_g + \spanR \{ v_{\beta(1)},\cdots,v_{\beta(n-h)} \} = \gX$. The proof of Theorem~\ref{th:kchoice} can now be extended to this case.
\endproof

Condition (\ref{eq:condLh}) guarantees the existence of two sets of vectors $\{v_1,\ldots,v_{n-h}\}$ and $\{w_1,\ldots,w_{n-h}\}$ {such that $\gV^\star_g + \spanR \{ v_{1},\cdots,v_{n-h}\} = \gX$ and $\bsmat A-\lambda_{\beta(i)}\,I_n & B \\[1mm] C_{({\beta(i)})} & D_{({\beta(i)})} \esmat \bsmat v_{i} \\[1mm] w_{i} \esmat=0$ for all $1\leq j\leq n-h$, where $\beta: \{1,\ldots,n-h\} \longrightarrow \delta$ is a bijective mapping. In such a case, $\rank [\,v_1\;\;\ldots\;\;v_{n-h}\;\;V_g\,]=n$ and we can compute $F=W\,V^{-1}$
which gives $\sigma(A+B\,F)=\{\lambda_{\beta(1)},\ldots,\lambda_{\beta(n-h)},\mu_1,\ldots,\mu_h\}$. The} tracking error is made up of at most a single closed-loop mode per component, but $h-(n-p)$ components of the tracking error are identically equal to zero, which means that in those components the output is identically equal to the corresponding component of the reference signal for any initial condition (and this obviously can only happen whenever the corresponding row of the feedthrough matrix $D$ is non-zero).
Although fully tractable, the necessary and sufficient condition proposed in Proposition~\ref{th:kchoiceGen} requires to test each $\binom{p}{n-h}$ possible injective map $\beta:\{1,\cdots,n-h\}\longrightarrow {\delta}$. 
A necessary and sufficient condition for Problem~\ref{prob:mono2} is
\begin{equation} 
\label{eq:condLhglob}
\dim\left(\gV_{g}^{\star}+\sum_{j\in{S}}\gR_{j}^{\star}(\lambda_{j})\right)\ge n-p+{\rm {card}({S})} \quad \forall\,{S}\in \{\mathfrak{S} \in 2^{\{1,\ldots,p\}}\,|\;\; \mathrm{card}\,{\mathfrak{S}}>h-(n-p)\},
\end{equation}
which clearly reduces to (\ref{eq:condL}) when $h=n-p$.
We omit the proof.

\section{Solution to Problem~\ref{prob:mono0}}

{In this section, the role played by the eigenvalues $\mathfrak{L}_p =(\lambda_1, \lambda_2, \cdots, \lambda_p)$ in the existence of solutions to Problem ~\ref{prob:mono2} is investigated.}

\begin{theorem}\label{th:Lchoice}
Problem~\ref{prob:mono0} admits solution if and only if 

\begin{equation}
\forall\,{S}\in2^{\{1,\ldots,p\}},\qquad\dim\left(\gV_{g}^{\star}+\sum_{j\in{S}}\gR_{j}^{\star}\right)\ge n-p+{\rm {card}({S})}. \label{eq:cond}
\end{equation}
\end{theorem}

\proof
Suppose that (\ref{eq:cond}) is not satisfied. This means that there exists ${S} \in 2^{\{1,\ldots,p\}}$ such that $\dim\left(\gV_{g}^{\star}+\sum_{j\in{S}}\gR_{j}^{\star}\right)<n-p+{\rm {card}({S})}$, which gives $\dim\left(\gV_{g}^{\star}+\sum_{j\in{S}}\gR_{j}^{\star}(\lambda_{j})\right)<n-p+{\rm {card}({S})}$ for any $( \lambda_{1},\cdots,\lambda_{p} )\in \hat{\mathfrak{L}}_p$, 
 since by \eqref{eq:gRdec} there holds $\gR_{j}^{\star}(\lambda_{j})\subseteq\gR_{j}^{\star}$ for all $j \in \{1,\cdots, p\}$ and $\lambda_{j} \in \mathbb{R}\backslash \gZ$. In view of Theorem~\ref{th:kchoice}, this shows that Problem~\ref{prob:mono2} is never solvable, which implies that Problem~\ref{prob:mono0} does not admit solution.

Let us now assume that (\ref{eq:cond}) is valid. Consider the $p$-tuples $(\lambda_{1},\cdots,\lambda_{p}) \in \hat{\mathfrak{L}}_{p}$ for which \eqref{eq:condL} does not hold, i.e., for which there exists ${S}\in2^{\{1,\ldots,p\}}$ satisfying
\begin{equation} 
\label{eq:dimLbad}
\mathfrak{q}: \qquad 	\dim\left(\gV_{g}^{\star}+\sum_{j\in{S}}\gR_{j}^{\star}(\lambda_{j})\right)<n-p+{\rm {card}({S})}.
\end{equation}
The set of all those $p$-tuples restricted to the subset $2^{\{1,\ldots,c\}} \subseteq 2^{\{1,\ldots,p\}}$, for $c\in \{1,\ldots,p\}$, is 
\begin{equation} \label{eq:Pcdef}
	\gP_{c} \defi \left\{ (\lambda_{1},\cdots,\lambda_{p})\in \hat{\mathfrak{L}}_{p} \,\Big|\,\exists{S}\in2^{\{1,\ldots,c\}}:\; \mathfrak{q} \right\}.
\end{equation}
We prove that $\gP_{p}$ has empty interior; indeed, in such case that there exists $( \lambda_{1},\cdots,\lambda_{p} )\in \hat{\mathfrak{L}}_p$ satisfying \eqref{eq:condL}, leading to a solution of Problem~\ref{prob:mono2} by virtue of Theorem~\ref{th:kchoice}. To prove this fact, we proceed by induction on $c\in \{1,\ldots,p\}$. Consider the following condition:
\begin{equation}
\forall\,{S}\in2^{\{1,\ldots,c\}},\qquad\dim\left(\gV_{g}^{\star}+\sum_{j\in{S}}\gR_{j}^{\star}\right)\ge n-p+{\rm {card}({S})}.\label{eq:condr}
\end{equation}
The Inductive Hypothesis (IH) for $c$ reads as
\begin{equation} \label{eq:IH}
\text{(IH)}:\; \eqref{eq:condr} \Rightarrow \gP_{c} \text{~has empty interior.}
\end{equation}

We show that (IH) holds for $c=1$,
i.e., if $\dim\left(\gV_{g}^{\star}+\gR_{1}^{\star}\right)\ge n-p+1$,
then $\gP_{1}$ has empty interior, where $\gP_{1}=\{(\lambda_{1},\cdots,\lambda_{p})\in \hat{\mathfrak{L}}_{p} \,|\dim\left(\gV_{g}^{\star}+\gR_{1}^{\star}(\lambda_{1})\right)<n-p+1\}$.
Suppose by contradiction that $\gP_{1}$ has non-empty interior. {Then, there exists an open interval contained in $\gP_{1}$ and hence there exists a set $L_{1}\subset \gP_{1}$ composed of $r_{1}=\dim \gR_1^\star$ distinct real numbers not coincident with the invariant zeros of $\Sigma$. By Assumption~\ref{assVg} and the definition of $\gP_{1}$, for all $i\in\{1,\ldots,r_{1}\}$, 
\begin{equation}
n-p \leq \dim (\gV_{g}^{\star})
 \leq \dim \left( \gV_{g}^{\star}+\gR_{1}^{\star}(\lambda_{1}^{i}) \right)
	 < n-p+1,
\end{equation}
where $\{ \lambda_1^1,\ldots ,\lambda_1^{r_1}\}=L_{1}$.
This} implies that $\dim (\gV_{g}^{\star})=\dim \left( \gV_{g}^{\star}+\gR_{1}^{\star}(\lambda_{1}^{i}) \right)=n-p$, and hence $\gR_{1}^{\star}(\lambda_{1}^{i})\subseteq\gV_{g}^{\star}$. Thus, $\gR_{1}^{\star}(\lambda_{1}^{1})+\ldots+\gR_{1}^{\star}(\lambda_{1}^{r_{1}})\subseteq\gV_{g}^{\star}$.
Since the elements of $L_{1}$ are distinct from the invariant zeros of $\Sigma$, \eqref{eq:gRdec} ensures that $\gR_{1}^{\star}=\gR_{1}^{\star}(\lambda_{1}^{1})+\ldots+\gR_{1}^{\star}(\lambda_{1}^{r_{1}})$.
This gives $\gR_{1}^{\star}\subseteq\gV_{g}^{\star}$, which in turn leads to $\dim\left(\gV_{g}^{\star}+\gR_{1}^{\star}\right)= \dim\left(\gV_{g}^{\star}\right) = n-p$. Since \eqref{eq:cond} immediately leads to
$\dim\left(\gV_{g}^{\star}+\gR_{1}^{\star}\right)\ge n-p+1$, we get to a contradiction. We conclude that $\gP_{1}$ has empty interior and (IH) is verified for $c =1$.
Next, let $c \in \{1, \dots, p-1\}$ and assume that (IH) holds for $c$; we show that (IH) also holds for $c+1$.
To this end, let us introduce 
\begin{equation} \label{eq:Hcdef}
	\gH_{c+1} \defi \left\{ (\lambda_{1},\cdots,\lambda_{p})\in \hat{\mathfrak{L}}_{p} \,\Big|\,\exists{S}\in2^{\{1,\ldots,c+1\}} \setminus \{1,\ldots,c+1\}:\;\mathfrak{q} \right\} \subset \gP_{c+1}.
\end{equation}
Observe that $\gP_{c+1}$ can be decomposed as $(\gP_{c+1}\cap\gH_{c+1}) \cup ( \gP_{c+1}\backslash\gH_{c+1} )$. Thus, to prove that $\gP_{c+1}$ has empty interior, it suffices to prove that both $\gP_{c+1}\cap\gH_{c+1}$ and $\gP_{c+1}\backslash\gH_{c+1}$ have empty interior. Corollary~\ref{cor:grid} {in Appendix~B} ensures that this is true for the latter. To prove that this also holds for the former, we first show that $\gH_{c+1}$ has empty interior. 
Let
\[
\gH_{c+1}^{(j)}\defi\left\{ (\lambda_{1},\cdots,\lambda_{p})\in \hat{\mathfrak{L}}_{p} \,\Big|\,\exists{S}\in2^{\{1,\cdots,c+1\}\setminus\{j\}}:\;\mathfrak{q}\right\}, 
\]
and consider the condition
\begin{equation}
\forall{S}\in2^{\{1,\cdots,c+1\}\setminus\{j\}},\quad
\dim\left(\gV_{g}^{\star}+\sum_{j \in{S}}\gR_{j}^{\star}\right) \geq h+{\rm {card}({S})}\label{eq:condLrp1-1bis}
\end{equation}
for all $j\in\{1,\cdots,c+1\}$. 
By means of a simple reindexing, it is seen that (IH) is equivalent to
\begin{equation} \label{eq:IHgen} 
	\eqref{eq:condLrp1-1bis} \Rightarrow \gH_{c+1}^{(j)} \text{~has empty interior},
\end{equation}
which is now valid for all $j\in\{1,\ldots,c+1\}$.
From the trivial identities
\begin{equation}
2^{\{1,\ldots,c+1\}}\setminus\{1,\cdots,c+1\} = \bigcup_{j=1}^{c+1}\left\{ {S}\in2^{\{1,\ldots,c+1\}}\,\Big|\, j\notin{S}\right\}
	= \bigcup_{j=1}^{c+1}2^{\{1,\cdots,c+1\}\setminus\{j\}}, \label{eq:2pw}
\end{equation}
we see that $\gH_{c+1}$ can be written as
\begin{equation} 
\label{eq:Hredef}
\gH_{c+1}=\left\{ (\lambda_{1},\cdots,\lambda_{p})\in \hat{\mathfrak{L}}_{p} \,\Big|\,\exists j\in\{1,\cdots,c+1\},\;\exists{S}\in2^{\{1,\cdots,c+1\}\setminus\{j\}}:\;\mathfrak{q} \right\},
\end{equation}
which leads to the decomposition
\begin{equation} \label{eq:Hcp}
\gH_{c+1}=\bigcup_{j=1}^{c+1} \gH_{c+1}^{(j)}.
\end{equation}
In view of \eqref{eq:2pw}, if the condition
\begin{equation} \label{eq:almostcond}
\forall{S}\in2^{\{1,\ldots,c+1\}}\setminus\{1,\cdots,c+1\},\quad \dim\left(\gV_{g}^{\star}+\sum_{j\in{S}}\gR_{j}^{\star}\right) \geq h+{\rm {card}({S})}
\end{equation}
is satisfied, then \eqref{eq:condLrp1-1bis} is satisfied for all $j=\{1,\cdots,c+1\}$.

The proof that $\gH_{c+1}$ has empty interior can now be established. First observe that \eqref{eq:cond} implies \eqref{eq:almostcond} and hence \eqref{eq:condLrp1-1bis}, for $j\in\{1,\ldots,c+1\}$. Second, \eqref{eq:IHgen} implies that for all $j\in\{1,\ldots,c+1\}$ the set $\gH^{(j)}_{c+1}$ has empty interior, which by \eqref{eq:Hcp} ensures that $\gH_{c+1}$ has empty interior as well. Thus, (IH) is verified for $c +1$. For the arbitrariness of $c$, \eqref{eq:IH} holds for $c = p$, i.e., if (\ref{eq:cond}) is satisfied, then $\gP_{p}$ has empty interior.
\endproof

As for Theorem \ref{th:kchoice}, it is not difficult at this point to see that when $\dim \gV^\star_g>n-p$, Problem~\ref{prob:mono0} admits solutions if and only if the condition
\begin{equation} 
\dim\left(\gV_{g}^{\star}+\sum_{j\in{S}}\gR_{j}^{\star}\right)\ge n-p+{\rm {card}({S})}
\end{equation}
holds true for all ${S}\in \{\mathfrak{S} \in 2^{\{1,\ldots,p\}}\,|\;\; \mathrm{card}\,{\mathfrak{S}}>h-(n-p)\}$.

%

\begin{remark}
Theorem~\ref{th:Lchoice} established that if (\ref{eq:cond}) is satisfied, the set of all $( \lambda_{1},\cdots,\lambda_{p} )\in \hat{\mathfrak{L}}_p$ for which \eqref{eq:condL} does not hold is {\em thin}, as it has empty interior. This is usually enough to guarantee that the elements of $\hat{\mathfrak{L}}_p$ for which Problem~\ref{prob:mono2} does not admit solution are, loosely speaking, pathological, since examples in which thin sets have non-zero Lebesgue measure have to be constructed {\em ad-hoc}, and can be considered as rarities. Nevertheless, at this stage it is only possible to
 conjecture that a stronger result holds, i.e., that the Lebesgue measure of this set within $\hat{\mathfrak{L}}_p$ - and hence within $\real^n$ - is zero. 
\end{remark}


\section*{Concluding remarks}
In this paper, the problem of achieving a monotonic step response from any initial condition has been addressed for the first time in the literature for LTI MIMO systems.
 This new approach opens the door to a range of developments that for the sake of conciseness cannot be addressed in this paper, but that we briefly discuss:
 \begin{itemize}
 \item In the case that global monotonicity cannot be achieved, it is important to find structural conditions ensuring that every component of the tracking error consists of the sum of at most two, three, or more closed-loop modes. In such case, even if the response is not globally monotonic, it is still monotonic starting from suitable initial conditions. Thus, an important issue is the characterisation of the regions of the state space where the initial state must belong to guarantee that the system response can be made monotonic;
 \item A second relevant problem
 involves the use of the method in \cite{CDC} and \cite{SICON} to the end of computing the state feedback that achieves a globally monotonic step response and which at the same time delivers a robust closed-loop eigenstructure, by ensuring that the closed-loop eigenvalues are rendered insensitive to perturbations in the state matrices. 
This task can be accomplished by
 obtaining a feedback matrix that minimises the Frobenius condition number of the matrix of closed-loop eigenvectors, which is a commonly used robustness measure. 
 The problem of obtaining a feedback matrix with minimum gain measure can be handled in a similar way, by minimising 
 the Frobenius norm of the feedback matrix. 
 
 \item Another important extension is the one of time-varying references, along the same lines of the extension that was proposed in \cite{Schmid-NG-CDC} of the procedure introduced in \cite{Schmid-N-09}.
 
\item Using the same approach of \cite{Schmid-N-09-2}, this method can be extended to the case of multivariable dynamic
output feedback tracking controllers.
\end{itemize}

\section*{Appendix A}

Let $h$ be the dimension of $\gV_{g}^{\star}$. Throughout this Appendix, we consider an arbitrary integer $k\in\{1,\ldots,n-h\}$ and a set of $k$ non-zero subspaces of $\gX$ denoted by $\{\gM_{1},\cdots,\gM_{k}\}$.

\begin{definition}
Let us introduce the following proposition:
\begin{equation}
\mathfrak{p}_1^k: \qquad \dim\left( \gV_{g}^{\star}+\spanR\{v_{1},\cdots,v_{k}\}\right) <h+k.\label{eq:condPath}
\end{equation}
The proposition $\neg\mathfrak{p}_1^k$ 
corresponds to the condition $\dim ( \gV_{g}^{\star}+\spanR\{v_{1},\cdots,v_{k}\} ) \ge h+k$. Since clearly $h+k$ vectors cannot span a subspace of dimension strictly greater than $h+k$, with a slight abuse for the sake of simplicity we will consider that 
\begin{equation*}
\neg \mathfrak{p}_1^k : \qquad \dim\left( \gV_{g}^{\star}+\spanR\{v_{1},\cdots,v_{k}\}\right) = h+k.
\end{equation*}
Using those equations, we define the sets
\begin{eqnarray}
\gK_{k}&\defi&\left\{ (v_{1},\cdots,v_{k})\in\gM_{1}\times\cdots\times\gM_{k}\big|\mathfrak{p}_1^k\right\}; \label{eq:path1} \\
\gK_{k}^c&\defi&\left\{ (v_{1},\cdots,v_{k})\in\gM_{1}\times\cdots\times\gM_{k}\big|\neg\mathfrak{p}_1^k\right\}. \label{eq:path2}
\end{eqnarray}
\end{definition}


\begin{lemma}
\label{lem:kchoice}
Let
\begin{equation}
\forall\,{S}\in2^{\{1,\ldots,k\}},\qquad\dim\left(\gV_{g}^{\star}+\sum_{j\in{S}}\gM_{j}\right)\ge h+{\rm {card}({S})}. \label{eq:condrlem}
\end{equation}
The following statements hold true:
\begin{enumerate}
\item If (\ref{eq:condrlem}) is satisfied, the set $\gK_{k}$ has measure zero;
\item If (\ref{eq:condrlem}) is not satisfied, then (\ref{eq:condPath})
holds for all $(v_{1},\cdots,v_{k})\in\gM_{1}\times\cdots\times\gM_{k}$.
\end{enumerate}
\end{lemma}

As a preliminary step toward the proof of Lemma~\ref{lem:kchoice}, consider the following definition. 

\begin{definition}
Let 
\begin{eqnarray}
\gA &\defi & \bigcap_{(v_{1},\cdots,v_{k}) \in\gK_{k}^c}\left(\gV_{g}^{\star}+\spanR\{v_{1},\cdots,v_{k}\}\right), \label{eq:Adef} \\
\gamma &\defi& \left\{ j \in \{1,\ldots,k\} \big|\gM_{j}\subseteq\gA \right\}, \label{eq:gammadef} \\
l &\defi& \mathrm{card}( \gamma ). \label{eq:ldef}
\end{eqnarray}
\end{definition}

\begin{proposition}
\label{fact:Lor2}
Given an arbitrary subspace $\gG\subseteq\gX$, define the following proposition:
\begin{equation}
\mathfrak{p}_2: \qquad \;\gG\subseteq\gV_{g}^{\star}+\spanR\{v_{1},\cdots,v_{k}\}.\label{eq:ninA}
\end{equation}
If the set $\{ (v_{1},\cdots,v_{k})\in\gK_{k}^c\big|\neg\mathfrak{p}_2 \}$ is not empty, then $\{ (v_{1},\cdots,v_{k})\in\gK_{k}^c\big|\mathfrak{p}_2\}$ has measure zero. 
\end{proposition}
This is a consequence of the fact that the Lebesgue measure of a proper subspace of a given a vector space is equal to zero.

\begin{proposition} \label{fact:Adim}
If $\gK_k^c$ is non-empty, then 
\begin{equation}
\dim\left( \gA\right) \leq h+l. \label{eq:ubAF}
\end{equation}
\end{proposition}

\proof
Let $\gamma^c \defi \left\{ j \in \{1,\ldots,k\} \big|\gM_{j}\nsubseteq\gA \right\}$. 
As a preliminary step, we prove that there exist $\{z_{j}\}_{j\in \{1,\ldots,k\}}$ and $\{ w_j\}_{j\in \gamma^{c}}$
such that 
\begin{eqnarray}
(z_{1},\cdots,z_{j-1},w_{j},z_{j+1},\cdots,z_{k})&\in&\gK_{k}^c \label{eq:zw1} \\
w_{j} &\notin& \gB_0 \label{eq:wj} \label{eq:zw2}
\end{eqnarray}
for all $j\in\gamma^{c}$, where
\begin{equation} \label{eq:B0def}
 \gB_0 \defi \gV_{g}^{\star}+\spanR\{z_{1},\cdots,z_{k}\}.
\end{equation}
To this end, define
\begin{eqnarray} 
\mathfrak{p}_3:\;\gM_{j} \ns&\ns\subseteq\ns&\ns \gV_{g}^{\star}+\spanR\{v_{1},\cdots,v_{k}\} \label{eq:Mjinc} \\
\gE_{k} \ns&\ns\defi\ns&\ns \left\{ (v_{1},\cdots,v_{k})\in\gK_{k}^c\big|\exists j\in\gamma^{c}:\; \mathfrak{p}_3 \right\} \\
\gE_{k}^{(j)} \ns&\ns\defi\ns&\ns \left\{ (v_{1},\cdots,v_{k})\in\gK_{k}^c\big| \mathfrak{p}_3 \right\} 
\end{eqnarray}
and let us first prove that {${m}_{\scriptscriptstyle \gK_k^c}(\gE_{k})=0$}. By definition of $\gamma^c$ and $\gA$, observe that, if $\gK_k^c$ is non-empty then the set $\{(v_{1},\cdots,v_{k})\in\gK_{k}^c | \neg \mathfrak{p}_3\}$ is non-empty for all $j\in\gamma^{c}$. In such a case, Lemma~\ref{fact:Lor2} ensures that ${m}_{\scriptscriptstyle \gK_k^c}(\gE_{k}^{(j)})=0$ which leads to ${m}_{\scriptscriptstyle \gK_k^c}(\gE_{k})=0$ by ${m}_{\scriptscriptstyle \gK_k^c}(\gE_{k})\leq\sum_{j\in\gamma^{c}}{m}_{\scriptscriptstyle \gK_k^c}(\gE_{k}^{(j)})=0$, which in turn follows from $\gE_{k}=\bigcup_{j\in\gamma^{c}}\gE_{k}^{(j)}$. This guarantees that there exists a particular element of $\gK_k^c$, denoted by $(z_{1},\cdots,z_{k})$, which does not belong to $\gE_{k} \subseteq \gK_k^c$. Hence, $\gM_{j}\nsubseteq\gV_{g}^{\star}+\spanR\{z_{1},\cdots,z_{k}\}=\gB_{0}$
 for all $j\in\gamma^{c}$ by definition of $\gE_{k}$. It readily follows that for all $j\in\gamma^{c}$ there exists $w_{j}\in\gM_{j}$ satisfying \eqref{eq:zw2}. It remains to prove that \eqref{eq:zw1} is verified for all $j\in\gamma^{c}$. This follows by observing that \textbf{(i)} by construction $(z_{1},\cdots,z_{j-1},w_{j},z_{j+1},\cdots,z_{k})\in\gM_{1}\times\cdots\times\gM_{k}$
and \textbf{(ii)} $\dim\left( \gV_{g}^{\star}+\spanR\{z_{1},\cdots,z_{j-1},w_{j},z_{j+1},\cdots,z_{k}\}\right) =h+k$
since $\dim\left( \gV_{g}^{\star}+\spanR\{z_{1},\cdots,z_{k},w_{j}\}\right) =h+k+1$
because $(z_{1},\cdots,z_{k})\in\gK_k^c$ and $w_{j}\notin\gB_{0}$.

Now that the existence of vectors satisfying \eqref{eq:zw1} and \eqref{eq:zw2} has been established, we define 
\[
\gB_{j} \defi \gV_{g}^{\star}+\spanR\{z_{1},\cdots,z_{j-1},w_{j},z_{j+1},\cdots,z_{k}\}, \qquad j\in\gamma^{c}.
\]
Observe that $\gA\subseteq\gB_{0}\cap (\bigcap_{j\in\gamma^{c}}\gB_{j} )$. In fact, $\gA=\gA\cap\gB_{0}\cap(\bigcap_{j\in\gamma^{c}}\gB_{j} )\subseteq\gB_{0}\cap (\bigcap_{j\in\gamma^{c}}\gB_{j} )$
can be deduced from the definition of $\gA$ bearing in mind that $(z_{1},\cdots,z_{k})\in\gK_k^c$ and (\ref{eq:zw1}) holds for all $j\in\gamma^{c}$. 
In the following, we prove that
\begin{equation} \label{eq:inters}
\dim \left(\gB_{0}\cap (\bigcap_{j\in\gamma^{c}}\gB_{j} ) \right) = h+l,
\end{equation}
which, in turn, implies \eqref{eq:ubAF}. 
To this end, we first show that for all $(\delta,j)$ such that $\delta\subseteq\{1,\cdots,k\}$, $j\in\gamma^{c}\cap \delta$, we have
\begin{equation} \label{eq:lemdelta}
(\gV_{g}^{\star}+\sum_{i\in\delta}\spanR\{z_{i}\})\cap\gB_{j}=
\gV_{g}^{\star}+\sum_{i\in\delta\backslash\{j\}}\spanR\{z_{i}\}.
\end{equation}
Define $\gH\defi\gV_{g}^{\star}+\sum_{i\in\delta}\spanR\{z_{i}\}$. 
Observe that $\gV_{g}^{\star}+\sum_{i\in\delta\backslash\{j\}}\spanR\{z_{i}\}$ is a subspace of both $\gH$ and $\gB_{j}$, so it is contained in their intersection. 
Moreover, $\gV_{g}^{\star}+\sum_{i\in\delta\backslash\{j\}}\spanR\{z_{i}\}$ and $\gH\cap \gB_{j}$ have the same dimension, which gives \eqref{eq:lemdelta}.
Indeed, $\dim(\gV_{g}^{\star}+\sum_{i\in\delta\backslash\{j\}}\spanR\{z_{i}\})=h+{\mathrm{card}}(\delta)-1$ since $\dim (\gB_0)=h+k$.
Using the Grassman rule, we have $\dim(\gH\cap\gB_{j})=\dim\gH+\dim\gB_{j}-\dim(\gH+\gB_{j})$, which reduces to
$\dim(\gH\cap\gB_{j})=(h+{\mathrm{card}}(\delta))+(h+k)-(h+k+1)=h+{\mathrm{card}}(\delta)-1$ because 
$\gH+\gB_{j}=\gB_{0}+\spanR\{w_{j}\}$ since $j\in\delta$ and $\dim(\gB_{0}+\spanR\{w_{j}\})=h+k+1$
since $w_{j}\notin\gB_{0}$. Then, applying \eqref{eq:lemdelta} with $\delta =\{1,\cdots,k\}$, we for all $j_1\in \gamma^c$ there holds $\gB_{0}\cap\gB_{j_1}=\gV_{g}^{\star}+\sum_{i\in\{1,\cdots,k\}\backslash\{j_1\}}\spanR\{z_{i}\}$. 
Similarly, it can be established that for all $\{j_1,j_2\}\subseteq \gamma^c$ we have $\left(\gB_{0}\cap\gB_{j_1}\right)\cap\gB_{j_2}=\gV_{g}^{\star}+\sum_{i\in\{1,\cdots,k\}\backslash\{j_1,j_2\}}\spanR\{z_{i}\}$
by \eqref{eq:lemdelta} with $\delta =\{1,\cdots,k\}\backslash\{j_1\}$. By repeating the same procedure, we obtain
\begin{equation*}
\gB_{0}\cap\left( \bigcap_{j\in\gamma^{c}}\gB_{j} \right)=\gV_{g}^{\star}+\sum_{i\in\{1,\cdots,k\}\backslash \gamma^{c} }\spanR\{z_{i}\}=\gV_{g}^{\star}+\sum_{i\in\gamma}\spanR\{z_{i}\}.
\end{equation*}
Then, \eqref{eq:inters} - and hence \eqref{eq:ubAF} - follows readily by observing that $\dim( \gB_0) =h+k$ implies $\dim( \gV_{g}^{\star}+\sum_{j\in\gamma}\spanR\{z_{i}\} ) =h+l$.
\endproof

In light of Proposition~\ref{fact:Lor2} and Proposition~\ref{fact:Adim}, the proof of Lemma~\ref{lem:kchoice} can now be established.

\begin{proofof}{\em Lemma~\ref{lem:kchoice}:}
We first prove the second point. Suppose that (\ref{eq:condrlem}) is not satisfied, i.e., there exists ${S}\in2^{\{1,\ldots,k\}}$ such that $\dim(\gV_{g}^{\star}+\sum_{j\in{S}}\gM_{j}) < h+\mathrm{card}({S})$. This implies that, for every collection of vectors $v_j$ such that $v_{j}\in \gM_j$ for all $j\in{S}$, we have $\dim(\gV_{g}^{\star}+\sum_{j\in{S}}\spanR\{v_{j}\})< h+\mathrm{card}({S})$ 
since $\gV_{g}^{\star}+\sum_{j\in{S}}\spanR\{v_{j}\}\subseteq \gV_{g}^{\star}+\sum_{j\in{S}}\gM_{j}$. This means that there exists a linear dependence among vectors $v_j$ and any basis of $\gV_{g}^{\star}$. Consequently, $\mathfrak{p}_1$ is satisfied for all $(v_{1},\cdots,v_{k})\in\gM_{1}\times\cdots\times\gM_{k}$.

We now assum that (\ref{eq:condrlem}) holds. For brevity, let $\tilde{\gM}_j \defi \gM_{1}\times\cdots\times\gM_{j}$. Let us prove that { ${m}_{\scriptscriptstyle \tilde{\gM}_k}(\gK_{k})=0$} by induction on $i\in \{1,\ldots,k\}$. The Inductive Hypothesis (IH) for $i$ reads: { ${m}_{\scriptscriptstyle \tilde{\gM}_i}(\gK_{i})=0$}. We first prove (IH) for $i=1$. Observe that \eqref{eq:condrlem} implies $\dim ( \gV_{g}^{\star}+\gM_{1}) \geq h+1$ and hence $\gM_{1}\nsubseteq\gV_{g}^{\star}$. Consequently, a generic vector $v_{1}\in\gM_{1}$ satisfies $v_{1}\notin\gV_{g}^{\star}$. This is equivalent to saying that { ${m}_{\scriptscriptstyle \gM_{1}}(\gK_{1})=0$}.
Let $i \in \{1, \dots, k-1\}$ and assume that { ${m}_{\scriptscriptstyle \tilde{\gM}_i}(\gK_{i})=0$}; we now show that { ${m}_{\scriptscriptstyle \tilde{\gM}_{i+1}}(\gK_{i+1})=0$}. Let us first introduce the following proposition
\begin{equation} \label{eq:condMp1}
\mathfrak{p}_4: \qquad \;\gM_{i+1}\subseteq\gV_{g}^{\star}+\spanR\{v_{1},\cdots,v_{i}\},
\end{equation}
and define 
\begin{eqnarray}
\gS_{i}&\defi&\left\{ (v_{1},\cdots,v_{i})\in\gM_{1}\times\cdots\times\gM_{i}\big|\mathfrak{p}_1^i\lor\mathfrak{p}_4\right\}, \\
\gS_{i}^c&\defi&\left\{ (v_{1},\cdots,v_{i})\in\gM_{1}\times\cdots\times\gM_{i}\big|\neg(\mathfrak{p}_1^i\lor\mathfrak{p}_4)
=\neg\mathfrak{p}_1^i\land\neg\mathfrak{p}_4
\right\}.
\end{eqnarray}
In the rest of the proof, we use the chain of implications
\begin{equation}
\eqref{eq:condrlem} \;\; \Rightarrow \;\; {m}_{\scriptscriptstyle \tilde{\gM}_i}(\gS_{i}\setminus\gK_{i})=0 \;\; \Rightarrow \;\; {m}_{\scriptscriptstyle \tilde{\gM}_i}(\gS_{i})=0 \;\; \Rightarrow \;\; {m}_{\scriptscriptstyle \tilde{\gM}_{i+1}}(\gK_{i+1})=0.
\end{equation}
To prove that ${m}_{\scriptscriptstyle \tilde{\gM}_i}(\gS_i\setminus\gK_i)=0$, it suffices to prove that $\gS^c_i$ is non-empty. This follows from Lemma~\ref{fact:Lor2} by observing that $\gS^c_i=\{ (v_{1},\cdots,v_{i})\in\gK_i^c\big|\neg\mathfrak{p}_4\}$ and $\gS_i\setminus\gK_i=\{ (v_{1},\cdots,v_{i})\in\gK_i^c\big|\mathfrak{p}_4\}$ because $(\mathfrak{p}_1^i\lor\mathfrak{p}_4)\land\neg\mathfrak{p}_1^i=\mathfrak{p}_4\land\neg\mathfrak{p}_1^i$. Suppose by contradiction that $\gS^c_i$ is empty.
First note that $\gK_i^c$ is non-empty because (IH) holds together with the decomposition $\gM_{1}\times\cdots\times\gM_{i}=\gK_i^c\cup\gK_i$, where every $\gM_j$ is a non-zero subspace.
Hence, the set $\gS^c_i$ -- where $\neg\mathfrak{p}_1^i\land\neg\mathfrak{p}_4$ holds -- is empty by assumption, whereas $\gK_i^c$ -- where only $\neg\mathfrak{p}_1^i$ holds -- is not. Thus, every element of $\gK_i^c$ satisfies (\ref{eq:condMp1}), which is equivalent to saying that $\gM_{i+1}\subseteq\gA$ where $\gA$ is given by \eqref{eq:Adef}.
Using this inclusion, observing that $\gV_{g}^{\star}\subseteq\gA$ and using the definition of $\gamma$ given by \eqref{eq:gammadef}, we get
\begin{equation}
\gV_{g}^{\star}+\sum_{j\in\gamma}\gM_{j}+\gM_{i+1}\subseteq\gA. \label{eq:AFfloor}
\end{equation}
By Lemma~\ref{fact:Adim}, the dimension of the subspace on the left hand-side of \eqref{eq:AFfloor} is smaller or equal to $h+l$. 
On the other hand, \eqref{eq:condrlem} ensures that -- for ${S}=\gamma \cup \{k+1\}$ -- this particular dimension is greater or equal to $h+l+1$, leading to a contradiction.
\footnote{
At this stage, we can even conclude that $\gA=\gV_{g}^{\star}+\sum_{j\in\gamma}\gM_{j}$.
Indeed, (\ref{eq:condrlem}) implies that $\dim(\gV_{g}^{\star}+\sum_{j\in\gamma}\gM_{j})\geq h+l$
which, together with $\gV_{g}^{\star}+\sum_{j\in\gamma}\gM_{j}\subseteq\gA$ - which is deduced from (\ref{eq:AFfloor}) - and (\ref{eq:ubAF}), allows
to write $h+l\leq\dim(\gV_{g}^{\star}+\sum_{j\in\gamma}\gM_{j})\leq\dim\left( \gA\right) \leq h+l$.
This leads to $\dim(\gV_{g}^{\star}+\sum_{j\in\gamma}\gM_{j})=\dim\left( \gA\right) =h+l$
and $\gA=\gV_{g}^{\star}+\sum_{j\in\gamma}\gM_{j}$.
}
Consequently, $\gS^c_i$ is non-empty and ${m}_{\scriptscriptstyle \tilde{\gM}_i}(\gS_i\setminus\gK_i)=0$.
The equality ${m}_{\scriptscriptstyle \tilde{\gM}_i}(\gS_{i})=0$ follows from \textbf{(i)} $\gS_i=(\gS_i\setminus\gK_i)\cup(\gS_i\cap\gK_i)$ which gives ${m}_{\scriptscriptstyle \tilde{\gM}_i}(\gS_i)\leq{m}_{\scriptscriptstyle \tilde{\gM}_i}(\gS_i\setminus\gK_i)+{m}_{\scriptscriptstyle \tilde{\gM}_i}(\gS_i\cap\gK_i)$, and \textbf{(ii)} ${m}_{\scriptscriptstyle \tilde{\gM}_i}(\gS_i\cap\gK_i)=0$. Indeed, ${m}_{\scriptscriptstyle \tilde{\gM}_i}(\gS_i\cap\gK_i)\leq{m}_{\scriptscriptstyle \tilde{\gM}_i}(\gK_i)=0$ as (IH) holds. Since ${m}_{\scriptscriptstyle \tilde{\gM}_i}(\gS_i)=0$, a generic $(v_{1},\cdots,v_i)\in\gM_{1}\times\cdots\times\gM_i=\gS_i \cup \gS_i^c$ belongs to $\gS_i^c$ and in turn verifies both $\neg\mathfrak{p}_1^i$ and $\neg\mathfrak{p}_4$. This easily gives ${m}_{\scriptscriptstyle \tilde{\gM}_{i+1}}(\gK_{i+1})=0$. Thus (IH) is valid for $i+1$ and hence for all $i \in \{1, \dots, k\}$.
\end{proofof}

\section*{Appendix B}
In this Appendix, we present a set of results that are used in the proof of Theorem~\ref{th:Lchoice}.

\begin{definition}
For any strictly positive integer $c$, let $L_1, \dots, L_c$ be finite non-empty sets of real numbers containing $l_1, \dots, l_c$ distinct elements, respectively. We say that { a set $G_c \subseteq \real^c$ is a {\em grid} in $\real^c$} of dimension $(l_1, \dots, l_c)$ if
$G_c = L_1 \times \dots \times L_c$. For each $j \in \{1, \dots, c\}$, we use
\be
L_j \defi \{ \lambda_j^{t_j}: 1 \leq t_j \leq l_j \} \label{indLp}
\ee
to denote an indexing of the elements of each $L_j$. We define $\Lambda = (\lambda_1, \dots, \lambda_c)$ as a node of the grid $G_c$, where each $\lambda_j = \lambda_j^{t_j}$, for some $1 \leq t_j \leq l_j$.
\end{definition}



\begin{lemma} \label{Badlambda}
For any $c\in \{1,\ldots,p-1\}$, let $L_1, \dots, L_{c+1}$ denote arbitrary sets containing $r_1, \dots, r_{c+1}$ distinct real numbers, respectively, not coincident with the invariant zeros of $\Sigma$.
Assume the sets $L_j$ are indexed as in (\ref{indLp}), and let $G_{c+1}$ 
be a grid in $\hat{\mathfrak{L}}_{p} \setminus \gH_{c+1}$ of dimension $(r_1, \dots, r_{p})$. Then, there exists (at least) one node $\Lambda \in G_{c+1}$ such that $\Lambda \notin \gP_{c+1}$.
\end{lemma}

\proof
Suppose by contradiction that every node $\Lambda$ of $G_{c+1}$ belongs to $\gP_{c+1}$, i.e., $G_{c+1} \subseteq \gP_{c+1}$. Since $G_{c+1} \subseteq \hat{\mathfrak{L}}_{c+1} \setminus \gH_{c+1}$, it follows that $G_{c+1} \subseteq \gP_{c+1} \cap (\hat{\mathfrak{L}}_{c+1} \setminus \gH_{c+1})$ which reduces to $G_{c+1} \subseteq \gP_{c+1} \setminus \gH_{c+1}$ as $\gP_{c+1} \subseteq \hat{\mathfrak{L}}_{c+1}$. From the definition of $\gP_{c+1}$ and $\gH_{c+1}$ given by \eqref{eq:Pcdef} by \eqref{eq:Hcdef}, respectively, for any $\Lambda\in G_{c+1}$, the only set ${S}\in2^{\{1,\ldots,c+1\}}$ for which \eqref{eq:dimLbad} holds is ${S}=\{1,\ldots,c+1\}$. This gives
\begin{equation} \label{eq:PCprop}
\dim \big(\gV_{g}^{\star}+\sum_{j\in{S}}\gR_{j}^{\star}(\lambda_{j}^{t_j})\big)
\text{~is~}\left\{ \begin{array}{ll}
\geq n-p+{\rm {card}({S})} & ({\rm {if~card}({S})} < c+1) \\
 < n-p+c+1 & ({\rm {if~card}({S})} = c+1) \\
\end{array} \right.
\end{equation}
for all $\Lambda \in G_{c+1}$.
Now, we define $K \defi \{1,\ldots,c+1\}$ and the subspace
\begin{equation}
\gW_c \defi \gV_{g}^{\star}+\sum_{j\in K} \ \gR_{j}^{\star}(\lambda_j^{1}).
\end{equation}
As an intermediate step, we want to prove that
\begin{equation} \label{eq:dimWc}
\dim( \gW_c) = n-p+c
\end{equation}
and
\begin{equation} \label{eq:Rloc}
\gR_{j}^{\star}(\lambda_j^{t_j}) \subseteq \gW_c 
\end{equation}
for all $j\in\{1,\cdots,c+1\}$ and for all $t_j \in \{1,\ldots,r_j\}$. To this end, let us first define 
\begin{equation}
\gW_c^{(j)} \defi \gV_{g}^{\star}+\sum_{j\in K\setminus \{j\}} \ \gR_{j}^{\star}(\lambda_j^{1}).
\end{equation}
Eq. \eqref{eq:PCprop} gives $\dim (\gW_c) < n-p+c+1$ and $\dim (\gW_c^{(j)}) \geq n-p+c$. 
Since clearly $\gW_c^{(j)} \subseteq \gW_c$, 
for all $j\in\{1,\cdots,c+1\}$ we find
\begin{equation}
n-p+c \leq \dim( \gW_c^{(j)})
	 \leq \dim (\gW_c)
	 < n-p+c+1.
\end{equation}
This leads to \eqref{eq:dimWc} and hence $\gW_c = \gW_c^{(j)}$ for all $j\in\{1,\cdots,c+1\}$. A similar argument can be used for the other nodes of $G_{c+1}$ {which can be expressed as $\Lambda=(\lambda_1^{1},\ldots,\lambda_{j-1}^{1},\lambda_{j}^{t_j}, \lambda_{j+1}^{1},\ldots,\lambda_{c+1}^{1})$ for some $1\leq j \leq c+1$ and $1\leq t_j \leq r_j$}. From $\gW_c^{(j)} + \gR_{j}^{\star}(\lambda_j^{t_j})=\gV_{g}^{\star}+\sum_{j\in{S}}\gR_{j}^{\star}(\lambda_{j}^{t_j})$ with ${S}=\{1,\ldots,c+1\}$, \eqref{eq:PCprop} gives
\begin{equation}
n-p+c \leq \dim (\gW_c^{(j)})
	 \leq \dim \left( \gW_c^{(j)} + \gR_{j}^{\star}(\lambda_j^{t_j}) \right)
	 < n-p+c+1
\end{equation}
for all $j\in\{1,\cdots,c+1\}$ and for all $t_j \in \{1,\ldots,r_j\}$ which clearly leads to \eqref{eq:Rloc} as $\gW_c = \gW_c^{(j)}$.

We now show that \eqref{eq:dimWc} and \eqref{eq:Rloc} contradict \eqref{eq:cond}. Since the elements of $L_{j} = \{ \lambda_{j}^{1}, \dots, \lambda_{j}^{r_{j}} \}$ are all distinct from the invariant zeros of $\Sigma$, we conclude from 
\eqref{eq:gRdec} and \eqref{eq:Rloc} that
$\gR_{j}^{\star}=\gR^{\star}_{j}(\lambda^{1}_{j})+\cdots+\gR^{\star}_{j}(\lambda^{r_{j}}_{j}) \subseteq \gW_c$. Because this hold for all $j\in \{1,\ldots,c+1\}$, it follows that
\begin{equation} \label{eq:incWc}
\gR_{1}^{\star} + \ldots +\gR_{c+1}^{\star} \subseteq \gW_c.
\end{equation}
On the other hand, from $\gR_{j}^{\star} (\lambda)\subseteq \gR_{j}^{\star}$ for all $\lambda \in \real \setminus \gZ$ we find $\gW_c \subseteq \gV^\star_g + \gR_{1}^{\star} + \ldots +\gR_{c+1}^{\star}$. This inclusion, together with the one obtained by adding $\gV^\star_g$ on both sides of \eqref{eq:incWc}, gives $\gW_c = \gV^\star_g + \gR_{1}^{\star} + \ldots +\gR_{c+1}^{\star}$.
According to \eqref{eq:dimWc}, the dimension of $\gW_c$ is $n-p+c$ which contradicts \eqref{eq:cond} ensuring that $\gW_c=\gV^\star_g + \sum_{j\in K} \gR_{j}^{\star} \geq n-p+c+1$. This allows to conclude that at least one node of $G_{c+1}$ does not belong to $\gP_{c+1}$.
\endproof


\begin{corollary} \label{cor:grid}
For any $c\in \{1,\ldots,p-1\}$, the set $\gP_{c+1}\setminus \gH_{c+1}$ has empty interior.
\end{corollary}

\proof
By contradiction, assume $\gP_{c+1}\setminus \gH_{c+1}$ has non-empty interior. Then there exists an open ball $U_{c+1} \subseteq \gP_{c+1}\setminus \gH_{c+1}$. This ball will contain a $c+1$-dimensional hypercube, and hence a $(r_1, \dots, r_{c+1})$-dimensional grid.
This contradicts Lemma \ref{Badlambda}, and hence $\gP_{c+1}\setminus \gH_{c+1}$ has empty interior.
\endproof

\end{document}